\documentclass[12pt]{article}
\usepackage{amssymb,amsmath,amsopn,amsfonts}
\usepackage[dvips]{color}
\usepackage{colordvi,multicol}
\usepackage{epsf}
\newfam\msbfam
\font\tenmsb=msbm10 \textfont\msbfam=\tenmsb \font\sevenmsb=msbm7
\scriptfont\msbfam=\sevenmsb \font\fivemsb=msbm5
\scriptscriptfont\msbfam=\fivemsb

\def\th#1{\vspace{1mm}\noindent{\bf #1}\quad}

\def\proof{\vspace{1mm}\noindent{\it Proof}\quad}

\numberwithin{equation}{section}

\def\bc{\begin{center}}
\def\ec{\end{center}}
\def\no{\noindent}
\def\hang{\hangindent\parindent}
\def\textindent#1{\indent\llap{\qquad #1\ \ \enspace}\ignorespaces}
\def\ref{\par\hang\textindent}

\begin{document}
\noindent{\LARGE\bf
Random attractor associated with\\\\  the quasi-geostrophic equation}
\vspace{0.5 true cm}

\noindent{\normalsize\sf
RongChan Zhu$^{1,3}$  XiangChan Zhu $^{2,3}$
\footnotetext{\baselineskip 10pt Research supported by  the DFG through IRTG 1132 and CRC 701\\
1.Department of Mathematics, Beijing Institute of Technology, Beijing 100081, China,\\
2.School of Sciences, Beijing Jiaotong University, Beijing 100044, China\\
3.Department of Mathematics, University of Bielefeld, D-33615 Bielefeld, Germany
}}

\vspace{0.2 true cm}

\noindent{
 Email:
 zhurongchan@126.com, z.xiangchan@gmail.com
\vspace{4mm}}%

\noindent {\bf Abstract}

We study the long time behavior of the solutions to the 2D stochastic quasi-geostrophic equation on $\mathbb{T}^2$ driven by additive noise and real linear multiplicative noise  in the subcritical case (i.e. $\alpha>\frac{1}{2}$) by proving the existence of a random attractor. The key point for the proof is the exponential decay of the $L^p$-norm and a boot-strapping argument.  The upper semicontinuity of random attractors is also established. Moreover, if the viscosity constant is large enough, the system has a trivial random attractor.

\vspace{1mm}
\no{\footnotesize{\bf 2010 Mathematics Subject Classification AMS}:\hspace{2mm} 60H15, 37L55,  35K55,  35Q86}
 \vspace{2mm}

\no{\footnotesize{\bf Keywords}:\hspace{2mm}   random attractors; quasi-geostrophic equation; random dynamical system; stochastic flow; stochastic partial differential equations.}

\section{Introduction}
Consider the following two dimensional (2D) stochastic quasi-geostrophic equation in the periodic domain $\mathbb{T}^2=\mathbb{R}^2/(2\pi \mathbb{Z})^2$:
$$\frac{\partial \theta(t,
\xi)}{\partial t}=-u(t,\xi)\cdot \nabla \theta(t,\xi)-\kappa (-\triangle)^\alpha \theta(t,\xi)+(G(\theta)\eta)(t,\xi),\eqno(1.1)$$
with  initial condition $$\theta(0,\xi)=\theta_0(\xi), \eqno(1.2)$$
where $\theta(t,\xi)$ is a real-valued function of $\xi\in \mathbb{T}^2$ and $t\geq0$, $0<\alpha<1, \kappa>0$ are real numbers. $u$ is determined by $\theta$ through a stream function $\psi$ via the following relations:
$$u=(u_1,u_2)=(-R_2\theta,R_1\theta)=R^\bot\theta.\eqno(1.3)$$

Here $R_j$ is the $j$-th periodic Riesz transform and $\eta(t,\xi)$ is a Gaussian random field, white noise in time, subject to the restrictions imposed below.
The case $\alpha=\frac{1}{2}$ is called the critical case, the case $\alpha>\frac{1}{2}$ sub-critical and the case $\alpha<\frac{1}{2}$ super-critical.

This equation is an important model in geophysical fluid dynamics. Indeed,
they are special cases of the general quasi-geostrophic approximations for atmospheric
and oceanic fluid flows with small Rossby and Ekman numbers. These models arise under
the assumptions of fast rotation, uniform stratification and uniform potential vorticity. The case $\alpha=1/2$ exhibits similar features (singularities) as the 3D Navier-Stokes equations and can therefore serve as a
 model case for the latter. In the deterministic case this equation has been intensively investigated because of both its mathematical importance and its background in geophysical fluid dynamics (see for instance \cite{CV10}, \cite{Re95}, \cite{Ju05}, \cite{KNV07} and the references therein).  In the deterministic case, the global existence of weak solutions has been obtained in \cite{Re95} and one most remarkable result in \cite{CV10} gives the existence of a classical solution for $\alpha=1/2$. In \cite{KNV07} another very important result is proved, namely that solutions for $\alpha=1/2$ with periodic $C^\infty$ data remain $C^\infty$ for all times. In the subcritical deterministic case, the long time behavior of the solution has been studied in \cite{Ju05} by proving the existence of the global attractor. In \cite{RZZ12} R${\ddot{o}}$ckner and the authors  studied the 2D stochastic quasi-geostrophic equation on
$\mathbb{T}^2$ for general parameter $\alpha\in (0,1)$ and for both additive as well as multiplicative noise case. For $\alpha>\frac{1}{2}$ R${\ddot{o}}$ckner and the authors obtained the existence and uniqueness of a (probabilistically strong) solution.

Recently there has been quite an interest in random attractors for stochastic partial differential equations. We refer the readers to \cite{BGLR11}, \cite{BL06}, \cite{CLR98}, \cite{C91} \cite{CDF97}, \cite{CF94}, \cite{FS96}, \cite{GLR11} and the references therein. Studying the global random attractor is one way to investigate the long time behavior of
partial differential equations perturbed by random noise. In this paper, we analyze the random attractor of the solutions
to  the stochastic quasi-geostrophic equation (1.1). More precisely,
we obtain that for the case $\alpha\in(\frac{1}{2},1)$, the random attractor  exists in the Soblev space $H^s$ (see definition below)  for any $s>2(1-\alpha)$ if the quasi-geostrophic equation is driven by additive noise (Theorem 3.7) or real linear multiplicative noise (Theorem 6.6). Moreover, the random attractor is infinitely smooth if the noise is sufficiently regular.

Comparing with some recent works on random attractors for SPDE (cf.\cite{BGLR11}, \cite{GLR11}), the main difficulty here lies in dealing with the nonlinear term in (1.1) since the dissipation term of the stochastic quasi-geostrophic equation is not  regular enough to control the nonlinear term as in the case of SPDE within the variational framework (see \cite{GLR11} for many examples).  In order to overcome this difficulty, we  consider the solution starting from a smaller state space $H^s$ (Sobolev space, see definition below) for $s>2(1-\alpha)$, which is an invariant subspace of the solution. We obtain a stochastic flow associated with the stochastic  quasi-geostrophic equation in $H^s$ space. Moreover, to get the existence of random attractors in $H^s$ space, one of the key point is
the improved positivity lemma we established in \cite{RZZ12} (see Lemma A.1). By this we obtain the decay of the $L^p$-norm of the solutions (cf. (3.16)), which is essential to obtain an absorbing ball in $H^s$.  On the other hand,  we can easily obtain  the exponential decay of the $L^2$-norm of the solution $\theta$ to the stochastic quasi-geostrophic equation, which implies an absorbing ball in $H$. However, this is not enough to obtain the existence of random attractor since this set is not compact in $H^s$ space. Here we apply the exponential decay of the integration $\int_t^{t+1}\|\theta(l)\|^2_{H^\alpha}dl$ and the decay of the $L^p$-norms of the solution to obtain  the exponential decay of the integration $\int_t^{t+1}\|\theta(l)\|^2_{H^{2\alpha}}dl$.  By using this and a similar technique we  obtain this kind of estimate for $\int_t^{t+1}\|\theta(l)\|^2_{H^{3\alpha}}dl$. Now we use a boot-strapping argument to conclude the existence of a compact absorbing ball in $H^s$ (Lemma 3.6).

Moreover, by the well-known results in \cite{CLR98} we  obtain the upper semi-continuity of the random attractors (Theorem 4.2) if the quasi-geostrophic equation is perturbed by a small $\varepsilon$-random perturbation, i.e. the random attractor is a random perturbation of the
deterministic one in the sense that, given a $\delta>0$, with probability one
there exists $\varepsilon_0$ (depending on $\omega$) sufficiently small, such that the random attractor is inside the
$\delta$-neighbourhood of the global attractor for all $\varepsilon<\varepsilon_0$.

 Furthermore, if the viscosity constant is large enough, we prove that the random attractor consists of a single point (Theorem 5.2). Since for the stochastic quasi-geostrophic equation the dissipation term is not strong enough to control the nonlinear term,
 we will use $L^p$-norm estimate to control the nonlinear term in a larger space.
  We first prove for almost every realization of the noise, trajectories starting from different initial conditions in $H^1$ converge to each other in a larger space $H^{-1/2}$ which is the dual space of Sobolev space $H^{1/2}$ (see Lemma 5.1). By this we obtain the existence of the limit for the stochastic flow $S(t,r,\omega)\theta_0$ constructed in Section 3 when time $r$ goes to $-\infty$. Then selecting a strictly stationarity version of the limiting process is the random attractor desired.

This paper is organized as follows. In Section 2 we recall some basic notions for random attractors and the stochastic quasi-geostrophic equation. In Sections 3 , we obtain the existence of a random attractor for the solutions of the stochastic quasi-geostrophic equation driven by  additive noise. In Section 4 we study the relation between the random attractor constructed in Section 3 and the global attractor obtained in \cite{Ju05} in the deterministic case, i.e. the upper semicontinuity of random attractors.  In Section 5 we obtain the system has a trivial random attractor if the viscosity constant is large enough. The existence of a random attractor for the solutions of the stochastic quasi-geostrophic equation driven by real multiplicative noise is established in Section 6.

\section{The basic set-up}
We first recall the notion of a random dynamical system (c.f.\cite{CF94}, \cite{CDF97}). Let $\{\vartheta_t:\Omega\rightarrow\Omega\}, t\in \mathbb{R}$, be a family of measure preserving transformations of a probability space $(\Omega,\mathcal{F},P)$ such that $(t,\omega)\mapsto\vartheta_t\omega$ is measurable, $\vartheta_0=id$ and $\vartheta_{t+r}=\vartheta_t\circ\vartheta_r$ for all $t,r\in\mathbb{R}$. Thus $((\Omega,\mathcal{F},P),(\vartheta_t)_{t\in\mathbb{R}})$ is a (measurable) dynamical system.
\vskip.10in
\th{Definition 2.1}  (i) A random dynamical system (RDS) on a Polish space $(X,d)$ with Borel $\sigma$-algebra $\mathcal{B}$ over $(\Omega,\mathcal{F},P,\vartheta_t)$ is a measurable map
$$\varphi:\mathbb{R}^+\times X\times \Omega\rightarrow X;(t,x,\omega)\mapsto\varphi(t,\omega)x$$
such that $\varphi(0,\omega)=id$ (identity on $X$) and
$$\varphi(t+r,\omega)=\varphi(t,\vartheta_r\omega)\circ\varphi(r,\omega),$$
for all $t,r\in \mathbb{R}^+$ and for all $\omega\in\Omega$. $\varphi$ is said to be a continuous RDS if $\varphi(t,\omega):X\rightarrow X$ is continuous for all $t\in \mathbb{R}^+$ and for all $\omega\in\Omega$.

(ii)A stochastic flow is a family of mappings $S(t,r,x;\omega):X\rightarrow X,-\infty<r\leq t<\infty$ parameterized by $\omega$ such that
$$(t,r,x,\omega) \rightarrow S(t,r;\omega)x$$
is $\mathcal{B}(\mathbb{R})\otimes\mathcal{B}(\mathbb{R})\otimes\mathcal{B}(X)\otimes\mathcal{F}/\mathcal{B}(X)$-measurable and
$$S(t,l;\omega)S(l,r;\omega)x = S(t,r;\omega)x,$$
$$S(t,r;\omega)x = S(t-r,0;\vartheta_r\omega)x,$$
for all $r\leq l\leq t$ and all $\omega\in \Omega$. $S$ is said to be a continuous stochastic flow if
$x\rightarrow S(t,r;\omega)x$ is continuous for all $r\leq t$ and $\omega\in\Omega$.
\vskip.10in

With the notion of an RDS above we can now recall the stochastic generalization of
notions of absorption, attraction and
$\Omega$-limit sets (cf. \cite{CF94}).
\vskip.10in

\th{Definition 2.2} (i) A (closed) set-valued map $K: \Omega\rightarrow 2^X$ is called measurable if $\omega\rightarrow K(\omega)$ takes values in the closed subsets of $X$ and for all $x\in X$ the map $\omega\mapsto d(x,K(\omega))$ is measurable, where for nonempty sets $A,B\in 2^X$ we set
$$d(A,B) = \sup\{\inf\{d(x,y):y\in B\},x\in A\}, d(x,B)=d(\{x\},B).$$
A measurable(closed) set-valued map is also called a (closed) random set.

(ii) Given a random set $K$, the set
$$\Omega(K,\omega)=\Omega_K(\omega)=\bigcap_{T\geq0}\overline{\bigcup_{t\geq T}\varphi(t,\vartheta_{-t}\omega)K(\vartheta_{-t}\omega)},$$
is said to be the $\Omega$-limit set of $K$.

(iii) Let $A, B$ be random sets. $A$ is said to absorb $B$ if $P$-a.s. there exists an absorption
time $t_B(\omega)$ such that for all $t\geq t_B(\omega)$
$$\varphi(t,\vartheta_{-t}\omega)B(\vartheta_{-t}\omega)\subset K(\omega).$$
$A$ is said to attract $B$ if $P$-a.s.
$$d(\varphi(t,\vartheta_{-t}\omega)B(\vartheta_{-t}\omega),A(\omega))\rightarrow0, \quad {t\rightarrow\infty}.$$

\vskip.10in
\th{Definition 2.3} A random attractor for an RDS  is a compact random set $A$ satisfying
$P$-a.s.:

(i) $A$ is invariant, i.e. $\varphi(t,\omega)A(\omega)=A(\vartheta_t\omega)$ for all $t > 0$.

(ii) $A$ attracts all deterministic bounded sets $B\subset X$.
\vskip.10in
The following proposition yields a sufficient criterion for the existence of a random attractor
of an RDS.
\vskip.10in
\th{Proposition 2.4}(cf. [8, Theorem 3.11]) Let $\varphi$ be an RDS on a Polish space $X$ and assume the existence of a compact random set $K$ absorbing every deterministic  bounded set $B\subset X$. Then there exists
a random attractor $A$, given by
$$A(\omega)=\overline{\bigcup_{B\subset X, B \textrm{ bounded }}\Omega_B(\omega)}.$$
\vskip.10in
In Section 3 and Section 6 we will apply Proposition 2.4 to prove the existence of a random attractor for the RDS
associated with the stochastic quasi-geostrophic equation.
\vskip.10in
Now we recall the following strong notion of stationarity, which is essential to construct an RDS.

\th{Definition 2.5} A map $Y:\mathbb{R}\times \Omega\rightarrow X$ is said to satisfy (crude) strict stationarity,  if
$$Y(t,\omega) = Y(0,\vartheta_t\omega),$$
for all $\omega\in\Omega$ and $t\in\mathbb{R}$ (for all $t\in \mathbb{R}$, $P$-a.s., where the zero-set may depend on $t$).
\vskip.10in

For RDS we need to use the following proposition  from [15,
Proposition 2.8] to select an indistinguishable
strictly stationary version.
\vskip.10in
\th{Proposition 2.6} Let $V\subset X$ and $Y:\mathbb{R}\times \Omega\rightarrow X$ be a process satisfying crude
strict stationarity. Assume that $Y\in C(\mathbb{R};X)\cap L^2_{\rm{loc}}(\mathbb{R};V)$ $P$-a.s.. Then there exists a process $\tilde{Y}:\mathbb{R}\times \Omega\rightarrow X$ such that

(i) $\tilde{Y}\in C(\mathbb{R};X)\cap L^2_{\rm{loc}}(\mathbb{R};V)$ for all $\omega\in\Omega$.

(ii) $Y, \tilde{Y}$ are indistinguishable, i.e.
$P[Y_t= \tilde{Y}_t, \textrm{ for any } t\in\mathbb{R}]=1$,
with a $\vartheta$-invariant exceptional set.

(iii) $\tilde{Y}$ is strictly stationary.
\vskip.10in

 In the following, we will restrict ourselves to flows which have zero average on the torus, i.e.
$$\int_{\mathbb{T}^2}\theta d\xi=0.$$
Thus (1.3) can be restated as
$$u=(-\frac{\partial \psi}{\partial \xi_2},\frac{\partial \psi}{\partial \xi_1}) \textrm { and } (-\triangle)^{1/2}\psi=-\theta.$$ Set $H=\{f\in L^2(\mathbb{T}^2):\int_{\mathbb{T}^2}f d\xi=0\}$ and let $|\cdot|$ and $\langle \cdot,\cdot\rangle$ denote the norm and inner product in $H$ respectively. On the periodic domain $\mathbb{T}^2$, $\{\sin (k\xi)|k\in \mathbb{Z}^2_+\}\cup\{\cos(k\xi)|k\in \mathbb{Z}^2_-\}$ form an eigenbasis of $-\triangle$ (we denote it by $\{e_k\}$). Here $ \mathbb{Z}^2_+=\{(k_1,k_2)\in \mathbb{Z}^2|k_2>0\}\cup\{(k_1,0)\in \mathbb{Z}^2|k_1>0\}, \mathbb{Z}^2_-=\{(k_1,k_2)\in \mathbb{Z}^2|-k\in \mathbb{Z}^2_+\}, \xi\in \mathbb{T}^2$, and the corresponding eigenvalues are $|k|^2$. Define $$\|f\|_{H^s}^2=\sum_k |k|^{2s}\langle f,e_k\rangle^2$$ and let $H^s$ denote the Sobolev space of all $f$ for which $\|f\|_{H^s}$ is finite. Set $\Lambda=(-\triangle)^{1/2}$. Then $$\|f\|_{H^s}=|\Lambda^s f|.$$

By the singular integral theory of Calder${\acute{o}}$n and Zygmund (cf [18, Chapter 3]), for any $l\geq0, p\in(1,\infty)$, there is a constant $C=C(l,p)$, such that
\begin{equation}\|\Lambda^lu\|_{L^p}\leq C(l,p)\|\Lambda^l\theta\|_{L^p}.\end{equation}

Fix $\alpha\in(0,1)$ and define the linear operator $A_\alpha: D(A_\alpha)=H^{2\alpha}(\mathbb{T}^2)\subset H\rightarrow H$ as $A_\alpha u:=\kappa (-\triangle)^\alpha u.$ The operator $A_\alpha$ is positive definite and self-adjoint with the same eigenbasis  as that of $-\bigtriangleup$ mentioned above. Denote the eigenvalues of $A_\alpha$ by $0<\lambda_1\leq\lambda_2\leq\cdots$ , and renumber the above eigenbasis correspondingly as $e_1,e_2,$....

First we recall the following important product estimates (cf. [16, Lemma A.4]):
\vskip.10in
\th{Lemma 2.7} Suppose that $s>0$ and $p\in (1,\infty)$. If $f,g\in C^\infty(\mathbb{T}^2)$, then
\begin{equation}\|\Lambda^s(fg)\|_{L^p}\leq C(\|f\|_{L^{p_1}}\|\Lambda^sg\|_{L^{p_2}}+\|g\|_{L^{p_3}}\|\Lambda^sf\|_{L^{p_4}}),\end{equation}
with $p_i\in (1,\infty), i=1,...,4$ such that
$$\frac{1}{p}=\frac{1}{p_1}+\frac{1}{p_2}=\frac{1}{p_3}+\frac{1}{p_4}.$$
\vskip.10in

We shall as well use the following Sobolev inequality (cf. [18, Chapter V]):
\vskip.10in
\th{Lemma 2.8} Suppose that $q>1, p\in [q,\infty)$ and
$$\frac{1}{p}+\frac{\sigma}{2}=\frac{1}{q}.$$
Suppose that $\Lambda^\sigma f\in L^q$, then $f\in L^p$ and there is a constant $C\geq 0$ such that
$$\|f\|_{L^p}\leq C\|\Lambda^\sigma f\|_{L^q}.$$

\th{Remark 2.9} Note that, because $div u=0$,  for regular functions $\theta$ and $\psi$, we have
$$\langle u(s)\cdot \nabla (\theta(s)+\psi),\theta(s)+\psi\rangle=0,$$ so
$$\langle u(s)\cdot \nabla \theta(s),\psi\rangle=-\langle u(s)\cdot \nabla \psi,\theta(s)\rangle.$$

\section{Additive noise}
 In this section we consider the abstract stochastic evolution equation driven by additive noise in place of Eqs (1.1)-(1.3),
\begin{equation}\frac{d\theta}{dt}+A_\alpha\theta+u\cdot\nabla \theta=dW,\end{equation}
where $u$ satisfies (1.3), $W$ is a trace-class two-sided Wiener process in $H$ with covariance $GG^*$ on a filtered probability space $(\Omega,\mathcal{F},\{\mathcal{F}_t\}_{t\in\mathbb{R}},P)$, where $G\in L_2(H,H)$ (i.e. $=$ all Hilbert-Schimit operators from $H$ to $H$.)

From now on we take $W$ to be the canonical process on $\Omega:=C_0(\mathbb{R},H):=\{w\in C(\mathbb{R},H);w(0)=0\}$, $\mathcal{F}_t$ to be canonical filtration and $\vartheta_t$ to be the Wiener shift given by $\vartheta_t\omega:=\omega(t+\cdot)-\omega(t)$ and $P=$ the law of $W$. Then $((\Omega,\mathcal{F},P),(\vartheta_t)_{t\in\mathbb{R}})$ is a (measurable) dynamical system.

In this section, we prove  that if the noise is regular, the associated random attractor is smooth. Now we
fix  $s>2(1-\alpha)$ and assume that:

\textbf{Hypothesis (E.1)} There exist $\varepsilon_0>0, \sigma_0>0\vee (1-s)$ such that $G\in L_2(H,H^{2+{\varepsilon_0}})\cap L_2(H,H^{s+1-\alpha+\sigma_0})$, i.e.
$$\mathcal{E}_0:=\rm{Tr}(\Lambda^{(2s+2-2\alpha+2\sigma_0)\vee (4+2\varepsilon_0)} GG^*)<\infty.$$

Given $\gamma>0$, let $z$ be the stationary solution of the equation:
$$dz+(A_\alpha+\gamma I)z=dW;$$
thus for $t\in\mathbb{R}$, $$z(t)=\int_{-\infty}^te^{-(t-l)(A_\alpha+\gamma I)} dW(l).$$

  By the Strong Law of Large Numbers ( see [9, Theorem 3.1.1]) and the assumption (E.1), we have for any $k\geq1$, $m\leq (s+1-\alpha+\sigma_0)\vee(2+\varepsilon_0)$
 \begin{equation}\lim_{t_0\rightarrow-\infty}\frac{1}{-1-t_0}\int_{t_0}^{-1}|\Lambda^{m}z|^kdl\rightarrow E|\Lambda^{m}z(0)|^k\quad P-a.s..\end{equation}
 Moreover,  $z\in C(\mathbb{R}, H^m)$ $P$-a.s.. By  Proposition 2.6 we can choose a version of $z$ such that it has strictly stationarity, i.e. for all $t\in \mathbb{R}$, $\omega\in \Omega$,
\begin{equation}z(t,\omega)=z(0,\vartheta_t\omega),\end{equation}
and for $\omega\in \Omega$, $z(\omega)\in C(\mathbb{R}, H^m)$. In the following we will take this version of  $z$. We can easily compute  \begin{equation}\lim_{\gamma\rightarrow\infty}E|\Lambda^{m}z(0)|^k=0,\end{equation} (cf. [3, Proposition 6.10]). By Ito's formula
 for $k\geq 2$ we have
$$\aligned &d|\Lambda^mz(t)|^k+k|\Lambda^mz(t)|^{k-2}|\Lambda^{m+\alpha} z|^2dt\\\leq& k|\Lambda^mz|^{k-2}\langle \Lambda^{2m}z,dW(t)\rangle+\frac{1}{2}k(k-1)|\Lambda^mz|^{k-2}\|\Lambda^mG\|^2_{L_2(H,H)}dt.\endaligned$$
By B-D-G inequality we can easily deduce that $m\leq(s+1-\alpha+\sigma_0)\vee(2+\varepsilon_0), k\geq 1$
$$E\sup_{0\leq t\leq 1}|\Lambda^mz(t)|^k\leq C(m,k).$$
Then by (3.3) and the dichotomy of linear growth (cf. [1, Proposition 4.1.3]) we have
\begin{equation}\limsup_{t\rightarrow\pm\infty}\frac{|\Lambda^mz(t)|^k}{|t|}=0,\end{equation}
on a $\vartheta$-invariant set of full $P$-measure.

\subsection{Stochastic flow}

 In the following we will consider the equation $\omega$-wise. If there is no confusion we omit $\omega$ for simplicity. We now use the change of variable $v(t)=\theta(t)-z(t)$. Then, formally, $v$ satisfies the equation
\begin{equation}\frac{dv}{dt}+A_\alpha v+u\cdot\nabla\theta=\gamma z.\end{equation}
For (3.6) we obtain the following $\omega$-wise existence and uniqueness result if the initial value starts from $H^s, s>2(1-\alpha)$.
\vskip.10in
\th{Theorem 3.1} Fix $\alpha>1/2$. Suppose the condition (E.1) holds. For any $v_0\in H^s, s>2(1-\alpha)$, there exists a unique solution $v\in L^\infty_{\rm{loc}}([t_0,\infty);H^s)\cap L^2_{\rm{loc}}([t_0,\infty); H^{s+\alpha})$ of equation (3.6) with $v(t_0)=v_0$, i.e. for any $\varphi\in C^1(\mathbb{T}^2)$
$$\langle v(t),\varphi\rangle-\langle v_0,\varphi\rangle+\int_{t_0}^t\langle A_\alpha^{1/2}v(r),A_\alpha^{1/2}\varphi\rangle dr-\int_{t_0}^t \langle (u_v+u_z)(r)\cdot \nabla \varphi,(v+z)(r)\rangle dr=\int_{t_0}^t\langle \gamma z,\varphi\rangle dr,$$
where $u_v,u_z$ satisfy (1.3) with $\theta$ replaced by $v,z$ respectively.

 \proof [Step 1] We first establish the existence and uniqueness of the solutions to the following linear equation:
\begin{equation}d v(t)+A_\alpha v(t)dt+w(t)\cdot \nabla (v(t)+z)dt= \gamma zdt\end{equation}
 $$v(t_0)\in H^\alpha\cap H^s,$$with a given smooth function $w(t)$ which satisfies $\rm{div} w(t)=0$ and $\sup_{t\in[t_0,T]}(\|w(t)\|_{C^{2}(\mathbb{T}^2)}+|\Lambda^{s+\alpha}w(t)|)\leq C(T), $ for any $T>t_0$.  Now consider the Galerkin approximation to (3.7):
 \begin{equation}dv^n(t)+A_\alpha v^n(t)dt+P_n(w(t)\cdot \nabla (v^n(t)+z))dt= P_n\gamma zdt,\end{equation}
 $$v^n(t_0)=P_nv_0,$$
 where $P_n$ is  the orthogonal projection in $H$ onto the linear space spanned by $e_1,...e_n$. Since all the coefficients are smooth in $P_nH$, this equation has a smooth solution $v^n$.
 We get the following estimate by taking the inner product in $L^2$ with $\Lambda^se_k$ for (3.8), multiplying both sides by $\langle v^n,\Lambda^se_k\rangle$ and summing up over $k$:
   $$\aligned\frac{1}{2}\frac{d}{dt}|\Lambda^sv^n|^2+\kappa|\Lambda^{s+\alpha} v^n|^2\leq&|\Lambda^{s-\alpha}(w\cdot\nabla(v^n+z))||\Lambda^{s+\alpha}v^n|+\gamma |\Lambda^s z||\Lambda^s v^n|\\\leq&C|\Lambda^{s+\alpha} v^n|[|\Lambda^{s-\alpha+1+\sigma_1}(v^n+z)|\|w\|_{L^{p_0}}+|\Lambda^{s-\alpha+1+\sigma_1}w|\|v^n+z\|_{L^{p_0}}]\\&+\gamma |\Lambda^s z||\Lambda^s v^n|\\\leq&C|\Lambda^{s+\alpha} v^n|^{1+r_0}|\Lambda^{s}v^n|^{1-r_0}\|w\|_{L^{p_0}}+C|\Lambda^{s+\alpha}v^n||\Lambda^{s-\alpha+1+\sigma_1}w|\|v^n\|_{L^{p_0}}\\&+\gamma |\Lambda^s z||\Lambda^s v^n|+C(T)|\Lambda^{s+\alpha}v^n|[|\Lambda^{s-\alpha+1+\sigma_1}z|+\|z\|_{L^{p_0}}]\\\leq&\frac{\kappa}{2}|\Lambda^{s+\alpha} v^n|^2+C(T)|\Lambda^sv^n|^2+C(T)\gamma^2|\Lambda^s z|^2+C(T)|\Lambda^{s-\alpha+1+\sigma_1}z|^2 ,\endaligned$$
    where $(1-s)\vee0<\sigma_1=2/p_0<(2\alpha-1)\wedge\sigma_0, r_0=\frac{1+\sigma_1-\alpha}{\alpha}$, $C(T)$ is a constant changing from line to line and we used Lemma 2.7 in the second inequality and the interpolation inequality in the third inequality and Young's inequality and $H^s\subset L^{p_0}$ in the last inequality.
By this estimate and $z\in C(\mathbb{R},H^{s-\alpha+1+\sigma_1})$, we get that $$\sup_{t\in[t_0,T]}|\Lambda^sv^n(t)|^2+\int_{t_0}^T|\Lambda^{s+\alpha}v^n(l)|^2dl\leq C,$$
where $C$ is a constant independent of $n$.
By a similar calculation and  $z\in C(\mathbb{R},H^{2})$ we also obtain that
$$\sup_{t\in[t_0,T]}|\Lambda^\alpha v^n(t)|^2+\int_{t_0}^T|\Lambda^{2\alpha}v^n(l)|^2dl\leq C,$$
where $C$ is a constant independent of $n$.
By (3.8) and the above estimates we know that
$$ \|v_n\|_{W^{1,2}([t_0,T],H)}\leq C.$$
By the compactness embedding $W^{1,2}([t_0,T],H^{-3})\cap L^2([t_0,T],H^{s+\alpha})\subset L^2([t_0,T],H^s)$ and $W^{1,2}([t_0,T],H)\subset C([t_0,T],H^{-1})$ we have that there exists a subsequence of $v^n$ converging in $L^2([t_0,T],H^s)\cap C([t_0,T],H^{-1})$ to a  function $v$ which is a solution to (3.7) and $v\in L^\infty([t_0,T];H^s\cap H^\alpha)\cap L^2([t_0,T],H^{s+\alpha}\cap H^{2\alpha})\cap C([t_0,T],H^{-1})$. Uniqueness of (3.7) is obvious.

[Step 2] We construct an approximation of (3.6) by a similar construction as in the proof of [17, Theorem 3.3]:

We pick a smooth $\phi\geq0$, with supp$\phi\subset[1,2], \int_0^\infty \phi=1$, and for $\delta>0$ let
$$U_\delta[\theta](t):=\int_0^\infty \phi(\tau)(k_\delta*R^\bot\theta)(t-\delta\tau)d\tau,$$ where $k_\delta$ is the periodic Poisson Kernel in $\mathbb{T}^2$ given by $\widehat{k_\delta}(\zeta)=e^{-\delta|\zeta|},\zeta\in \mathbb{Z}^2$, and we set $\theta(t)=0, t<t_0$. We take a zero sequence $\delta_n, n\in\mathbb{N}$, and consider the equation:
\begin{equation}\frac{d v_n(t)}{dt}+A_\alpha v_n(t)+u_n(t)\cdot \nabla (v_n(t)+z)= \gamma z\end{equation} with initial data $v_n(t_0)=k_{\delta_n}*v_0$ and $u_n=U_{\delta_n}[v_n+z]$.
For a fixed $n$, this is a linear equation in $v_n$ on each subinterval $[t_k,t_{k+1}]$ with $t_k=t_0+k\delta_n$, since $u_n$ is smooth and is determined by the values of $v_n$ on the two previous subintervals. By [Step 1], we obtain the existence of a  solution $v_n\in L^\infty([t_0,T];H^s\cap H^\alpha)\cap L^2([t_0,T],H^{s+\alpha}\cap H^{2\alpha})\cap C([0,T],H^{-1})$ to (3.9). Now for $s<1$, we choose $p$ such that $\frac{2}{(2\alpha-1)\wedge\sigma_0}<p\leq\frac{2}{1-s}$ and for $s\geq 1$ we take any $p$ satisfying $\frac{2}{(2\alpha-1)\wedge\sigma_0}<p<\infty$, where $\sigma_0$ appears in Assumption (E.1). From now on we fix such $p$ and we have $ H^s\subset L^p$ by Lemma 2.8. Since the periodic Riesz transform is bounded on $L^p$, we have for $t>t_0$ and $l\geq0$
\begin{equation}\sup_{[t_0,t]}\|\Lambda^lU_\delta[\theta]\|_{L^p}\leq C\sup_{[t_0,t]}\|\Lambda^l\theta\|_{L^p},\end{equation}
and also
\begin{equation}\int_{t_0}^t\|\Lambda^lU_\delta[\theta]\|_{L^p}^pd\tau\leq C\int_{t_0}^t\|\Lambda^l\theta\|_{L^p}^pd\tau.\end{equation}

By  Lemma A.1 we obtain for $v_n$ the following inequality by taking inner product with $|v_n|^{p-2}v_n$ in $L^2$

\begin{equation}\aligned\frac{d}{dt}\|v_n\|_{L^p}^p+2\lambda_1\|v_n\|_{L^p}^p\leq& p|\langle u_n\cdot\nabla (v_n+z) ,|v_n|^{p-2}v_n\rangle|+p\langle\gamma z,|v_n|^{p-2}v_n\rangle\\\leq &p\|\nabla z\|_{L^\infty}\|u_n\|_{L^p}\|v_n\|_{L^p}^{p-1}+Cp\gamma\|z\|_{L^p}\|v_n\|_{L^p}^{p-1},\endaligned\end{equation}
where we used $\rm{div}u_n=0$ and $\langle u_n\cdot\nabla v_n ,|v_n|^{p-2}v_n\rangle=0$ in the last inequality.
Therefore
$$\aligned&\|v_n(t)\|_{L^p}^p-\|v_n(t_0)\|_{L^p}^p+\int_{t_0}^t2\lambda_1 \| v_n(\tau)\|_{L^p}^pd\tau\\\leq &\varepsilon\int_{t_0}^t(\|u_n\|_{L^p}^p+\|v_n\|_{L^p}^p)d\tau+pC(\varepsilon)\int_{t_0}^t(\|\nabla z\|_{L^\infty}^{\frac{p}{p-1}}\|v_n\|_{L^p}^{p}+Cp\gamma\|z\|_{L^p}^p)d\tau\\\leq&\varepsilon\int_{t_0}^t\|v_n\|_{L^p}^pd\tau+pC(\varepsilon)\int_{t_0}^t(\|\nabla z\|_{L^\infty}^{\frac{p}{p-1}}\|v_n\|_{L^p}^{p}+Cp\gamma\|z\|_{L^p}^p)d\tau,\endaligned$$
where we used Young's inequality in the first inequality and (3.11) in the last inequality.
Then  Gronwall's  lemma, $\nabla z\in C(\mathbb{R},H^{1+\varepsilon})\subset C(\mathbb{R},L^\infty)$ for $\varepsilon\leq\varepsilon_0$ with $\varepsilon_0$ in (E.1) and $H^s\subset L^p$ yield that for any $T\geq t_0$
\begin{equation}\sup_{t\in[t_0,T]}\|v_n(t)\|_{L^p}\leq C,\end{equation}
where $C$ is a constant independent of $n$.

Moreover, we get the following estimate by taking the inner product in $L^2$ with $\Lambda^se_k$ for (3.9), multiplying both sides by $\langle v_n,\Lambda^se_k\rangle$ and summing up over $k$:
   $$\aligned\frac{1}{2}\frac{d}{dt}|\Lambda^sv_n|^2+\kappa|\Lambda^{s+\alpha} v_n|^2\leq&|\Lambda^{s-\alpha}(u_n\cdot\nabla(v_n+z))||\Lambda^{s+\alpha}v_n|+\gamma |\Lambda^s z||\Lambda^s v_n|\\\leq&C|\Lambda^{s+\alpha} v_n|[|\Lambda^{s-\alpha+1+\sigma_1}(v_n+z)|\|u_n\|_{L^p}+|\Lambda^{s-\alpha+1+\sigma_1}u_n|\|v_n+z\|_{L^p}]\\&+\gamma |\Lambda^s z||\Lambda^s v_n|,\endaligned$$
   where $\sigma_1=2/p<(2\alpha-1)\wedge\sigma_0$ and we used Lemma 2.7 in the last inequality.
Hence we obtain that for $r_0=\frac{1+\sigma_1-\alpha}{\alpha}, r=\frac{2\alpha}{2\alpha-1-\sigma_1}$,
\begin{equation}\aligned&\frac{1}{2}(|\Lambda^sv_n(t)|^2-|\Lambda^sv_n(t_0)|^2)+\kappa\int_{t_0}^t|\Lambda^{s+\alpha} v_n|^2d\tau\\\leq&C\int_{t_0}^t|\Lambda^{s+\alpha} v_n|[|\Lambda^{s-\alpha+1+\sigma_1}(v_n+z)|\|u_n\|_{L^p}+|\Lambda^{s-\alpha+1+\sigma_1}u_n|\|v_n+z\|_{L^p}]+\gamma |\Lambda^s z||\Lambda^s v_n|d\tau\\\leq&C\int_{t_0}^t[|\Lambda^{s+\alpha} v_n|^{1+r_0}|\Lambda^{s}v_n|^{1-r_0}+|\Lambda^{s+\alpha} v_n||\Lambda^{s-\alpha+1+\sigma_1}z|]\|u_n\|_{L^p}d\tau+\frac{\kappa}{4}\int_{t_0}^t|\Lambda^{s+\alpha} v_n|^2d\tau\\&+C\sup_{t\in[t_0,T]}\|v_n+z\|_{L^p}^2\int_{t_0}^t[|\Lambda^sv_n|^{2(1-r_0)}|\Lambda^{s+\alpha}v_n|^{2r_0}
+|\Lambda^{s-\alpha+1+\sigma_1}z|^2]d\tau+\int_{t_0}^t\gamma |\Lambda^s z||\Lambda^s v_n|d\tau\\
\leq&\frac{\kappa}{2}\int_{t_0}^t|\Lambda^{s+\alpha} v_n|^2d\tau+C[\sup_{t\in[t_0,T]}\|v_n+z\|_{L^p}^r+\|v_n+z\|^2_{L^p}+1]\int_{t_0}^t|\Lambda^sv_n|^2+
|\Lambda^{s-\alpha+1+\sigma_1}z|^2d\tau,\endaligned\end{equation}
where we used (3.10), (3.11), the interpolation inequality in the second inequality and Young's inequality in the last inequality.
By Gronwall's lemma, $z\in C(\mathbb{R};H^{s-\alpha+1+\sigma_1})$ and (3.13) we get that for $v_0\in H^s$
\begin{equation}|\Lambda^sv_n(t)|^2+\kappa\int_{t_0}^T|\Lambda^{s+\alpha} v_n|^2d\tau\leq C,\end{equation}
where $C$ is also a constant independent of $n$. By the same argument as above we obtain
$$\|v_n\|_{W^{1, 2}([t_0,T], H^{-3})}\leq C,$$
where $C$ is a constant independent of $n$.
 By the compactness embedding $W^{1,2}([t_0,T],H^{-3})\cap L^2([t_0,T],H^{s+\alpha})\subset L^2([t_0,T],H^s)$ we have that there exists a subsequence of $v_n$ converging in $L^2([t_0,T],H^s)$ to a  solution $v\in L^\infty_{\rm{loc}}([t_0,\infty);H^s)\cap L^2_{\rm{loc}}([t_0,\infty); H^{s+\alpha})$ of equation (3.6). Thus (3.15) is also satisfied for $v$. Uniqueness can be deduced from  a similar argument as in the proof [17, Theorem 5.1] (also see the proof of Theorem 3.3).$\hfill\Box$
\vskip.10in

Then  taking the limit for (3.12) and using Gronwall's lemma, we obtain the following estimate which is essential to get the existence of an absorbing set in $H^s$:
\begin{equation}\aligned&\|v(t)\|_{L^p}\leq \|v(t_0)\|_{L^p}\exp\{-\frac{2\lambda_1}{p}(t-t_0)+\int_{t_0}^t\|\nabla z(\tau)\|_{L^\infty} d\tau\}\\&+C\int_{t_0}^t(\|\nabla z(\tau)\|_{L^\infty}\|z(\tau)\|_{L^p}+C\gamma\|z(\tau)\|_{L^p})\exp\{-\frac{2\lambda_1}{p}(t-\tau)+\int_{\tau}^t\|\nabla z(l)\|_{L^\infty} dl\}d\tau, t\geq t_0.\endaligned\end{equation}

\vskip.10in
\th{Theorem 3.2} Fix $\alpha>1/2$. Suppose the condition (E.1) holds. The solution $v$ obtained in Theorem 3.1 is in $ C([t_0,\infty);H^s)$.

\proof Since $v\in L_{\rm{loc}}^2([t_0,\infty);H^{s+\alpha})$, by \cite{T84} it is sufficient to show that
$$\Lambda^s\frac{dv}{dt}\in L^2_{\rm{loc}}([t_0,\infty); H^{-\alpha}).$$
For $\varphi$ smooth enough, we have
$$\aligned|\langle\frac{dv}{dt},\Lambda^s \varphi\rangle|=&|\kappa\langle-\Lambda^{\alpha}v,\Lambda^{s+\alpha}\varphi\rangle-\langle (u\cdot \nabla(v+z)),\Lambda^s\varphi\rangle+\langle\gamma \Lambda^sz,\varphi\rangle|\\\leq&[\kappa|\Lambda^{s+\alpha}v|+C|\Lambda^{s-\alpha+1}(u\cdot(v+z))|]|\Lambda^\alpha \varphi|+\gamma|\Lambda^{s-\alpha}z||\Lambda^\alpha\varphi|\\\leq&C[|\Lambda^{s+\alpha}v|+|\Lambda^{s-\alpha+1+\sigma_1}
(v+z)|\|v+z\|_{L^p}+\gamma|\Lambda^{s-\alpha}z|]|\Lambda^\alpha \varphi|,\endaligned$$
where $(1-s)\vee0<\sigma_1=\frac{2}{p}<(2\alpha-1)\wedge\sigma_0$ as (3.14) and we used Lemma 2.7 in the last inequality.
Then by a similar calculation as (3.14)
$$\|\Lambda^s\frac{dv}{dt}\|_{H^{-\alpha}}\leq C(\|v+z\|_{L^p}+1)|\Lambda^{s+\alpha}v|+C\|v+z\|_{L^p}|\Lambda^{s-\alpha+1+\sigma_1}z|+\gamma|\Lambda^{s-\alpha}z|.$$
By (3.13), (3.15) $H^s\subset L^p$ and the regularity for $z$, we obtain for $-\infty<t_0<T<\infty$
$$\int_{t_0}^T\|\Lambda^s\frac{dv}{dt}(\tau)\|_{H^{-\alpha}}^2d\tau<\infty,$$
which implies that $v\in C([t_0,\infty);H^s)$.
$\hfill\Box$

\vskip.10in
\th{Theorem 3.3} Fix $\alpha>1/2$. Suppose the condition (E.1) holds. Then for any fixed $t>0, \omega\in \Omega$, the map $v_0\mapsto v(t,\omega,t_0, v_0)$ is  continuous  from $H^s$ into itself, where $v(t,\omega;t_0,v_0)$ is the solution of equation (3.6) with $v(t_0)=v_0$.

\proof Let $v_1,v_2$ be two solutions of (3.6) and $\zeta=v_1-v_2, \theta_1=v_1+z,\theta_2=v_2+z$. Then $\zeta$ satisfies the following equation:
$$(\frac{d}{dt}\zeta,\varphi)+\kappa(\Lambda^\alpha\zeta,\Lambda^\alpha \varphi)=-(u_1\cdot\nabla \zeta,\varphi)-(u_\zeta\cdot\nabla \theta_2,\varphi),$$
where $\varphi\in C^1(\mathbb{T}^2)$, $u_1,u_\zeta$ satisfy (1.3) with $\theta$ replaced by $\theta_1,\zeta$ respectively.

Taking $\varphi=\Lambda^{s}e_k$, multiplying both sides by $\langle \zeta,\Lambda^se_k\rangle$ and summing up over $k$  we have the following estimate since $v_i\in C([t_0,\infty);H^s)\cap L^2_{\rm{loc}}([t_0,\infty); H^{s+\alpha})$, $i=1,2$, by Theorems 3.1, 3.2
$$\aligned\frac{1}{2}\frac{d}{dt}|\Lambda^s\zeta|^2+\kappa|\Lambda^{s+\alpha} \zeta|^2=&-\langle\Lambda^s(u_1\cdot\nabla\zeta),\Lambda^s\zeta\rangle-\langle u_\zeta\cdot\nabla \theta_2,\Lambda^{2s}\zeta\rangle\\\leq&C|\Lambda^{s+\alpha} \zeta|[|\Lambda^{s-\alpha+1}(u_\zeta\theta_2)|+|\Lambda^{s-\alpha+1}(u_1\zeta)|]\\\leq&C|\Lambda^{s+\alpha} \zeta|[|\Lambda^{s-\alpha+1+\sigma_1}\zeta|\|\theta_2\|_{L^p}
+|\Lambda^{s-\alpha+1+\sigma_1}\theta_2|\|\zeta\|_{L^p}\\&+|\Lambda^{s-\alpha+1+\sigma_1}\theta_1|\|\zeta\|_{L^p}
+|\Lambda^{s-\alpha+1+\sigma_1}\zeta|\|\theta_1\|_{L^p}]\\\leq&C|\Lambda^{s+\alpha} \zeta|^{1+r_0}|\Lambda^{s}\zeta|^{1-r_0}[|\Lambda^s\theta_2|+|\Lambda^s\theta_1|]
\\&+|\Lambda^{s+\alpha} \zeta||\Lambda^s\zeta|[|\Lambda^{s-\alpha+1+\sigma_1}\theta_2|+|\Lambda^{s-\alpha+1+\sigma_1}\theta_1|]\\\
\leq&\frac{\kappa}{2}|\Lambda^{s+\alpha} \zeta|^2+C[|\Lambda^s\theta_2|^r+|\Lambda^s\theta_1|^r\\&+|\Lambda^{s+\alpha}v_2|^2+|\Lambda^{s-\alpha+1+\sigma_1}z|^2
+|\Lambda^{s+\alpha}v_1|^2]|\Lambda^s\zeta|^2,\endaligned$$
 where $r_0=\frac{1+\sigma_1-\alpha}{\alpha}, r=\frac{2\alpha}{2\alpha-1-\sigma_1}$ for some $(1-s)\vee0<\sigma_1=\frac{2}{p}<(2\alpha-1)\wedge\sigma_0$ as in (3.14)  and we used Lemmas 2.7 in the second inequality and Lemma 2.8, the interpolation inequality, $H^s\subset L^p$ in the third inequality and Young's inequality in the last inequality.
Then  Gronwall's lemma yields that
$$|\Lambda^s\zeta|^2\leq C|\Lambda^s\zeta(t_0)|^2\exp\{\int_{t_0}^T|\Lambda^s\theta_2(\tau)|^r+|\Lambda^s\theta_1(\tau)|^r+|\Lambda^{s-\alpha+1+\sigma}z|^2+|\Lambda^{s+\alpha}v_1(\tau)|^2+
|\Lambda^{s+\alpha}v_2(\tau)|^2d\tau\}.$$
Thus the result follows.$\hfill\Box$

\vskip.10in

Now for $\theta_0\in H^s$ we define
  $$\varphi(t,\omega)\theta_0:=v(t,\omega;0,\theta_0-z(0,\omega))+z(t,\omega), \quad t\geq0.$$
  $$S(t,r;\omega)\theta_0:=v(t,\omega;r,\theta_0-z(r,\omega))+z(t,\omega),\quad t,r\in \mathbb{R}.$$

  Combining Theorems 3.1-3.3 we obtain the following results.
  \vskip.10in
  \th{Theorem 3.4} Fix $\alpha>1/2$. Suppose the condition (E.1) holds. Then $\varphi(t,\omega)$ is a continuous random dynamical system and $S(t,r;\omega)$
is a continuous stochastic flow, which is called the stochastic flow associated with the stochastic quasi-geostrophic equation  driven by additive noise.

\proof By the $\omega$-wise uniqueness of the solution to equation
(3.6)  obtained in Theorem 3.1 and (3.3), we obtain that
$$S(t,r;\omega)=S(t,l;\omega)S(l,r;\omega),$$
$$S(t,r;\omega)x = S(t-r,0;\vartheta_r\omega)x,$$
$$\varphi(t+r,\omega)=\varphi(t,\vartheta_r\omega)\circ\varphi(r,\omega),$$
for all $t,l,r\in \mathbb{R}$ and for all $\omega\in\Omega$.
It remains to prove the measurability of $\varphi:\mathbb{R}\times \Omega\times H^s\rightarrow H^s$, which also implies the measurability of $S$ by the relation between $\varphi$ and $S$.  Since $\varphi(t,\omega)\theta_0=v(t,\omega;0,\theta_0-z(0,\omega))+z(t,\omega)$, $t\mapsto v(t,\omega;0,\theta_0)$ and $\theta_0\mapsto v(t,\omega;0,\theta_0)$
are continuous, we only need to prove the measurability of $\omega\mapsto v(t,\omega;0,\theta_0)$. By $\omega$-wise uniqueness of the solution to (3.6) we deduce that each subsequence of the convolution approximation $v_n(t,\omega;0,\theta_0)$ (which is measurable since $\omega$-wise uniqueness holds for (3.7)) we used in the proof of Theorem 3.1 has a subsequence  converging to the same $v(t,\omega;0,\theta_0)$ in $L^2([t_1,t_2],H^s)$ for some $t_1\leq t\leq t_2$. Thus we obtain that the whole sequence of $v_n(t,\omega;0,\theta_0)$  converges to $v(t,\omega;0,\theta_0)$ in $L^2([t_1,t_2],H^s)$, which implies the measurability of  $\omega\mapsto v(t,\omega;0,\theta_0)$.
$\hfill\Box$
\vskip.10in

\subsection{ Absorption in $H^s$ at time $t=-1$}
In this subsection we will prove the existence of an absorbing ball in the space $H^s$.
\vskip.10in
\th{Lemma 3.5} Suppose the condition (E.1) holds. There exists  random radius $r_1(\omega), c_1(\omega),c_2(\omega)>0$, such that for all $\rho>0$ there exists $t(\omega)\leq -1$ such that the following holds $P$-a.s.: For all $t_0\leq t(\omega)$ and all $\theta_0\in H^s$ with $|\Lambda^s \theta_0|\leq \rho$, the solution $v(t,\omega;t_0,\theta_0-z(t_0,\omega))$ with $v(t_0)=\theta_0-z(t_0,\omega)$ satisfies the following inequalities:
$$|\Lambda^sv(-1,\omega;t_0,\theta_0-z(t_0,\omega))|^2\leq r_1^2(\omega).$$
$$|\Lambda^sv(t,\omega;t_0,\theta_0-z(t_0,\omega))|^2\leq c_1(\omega), t\in[-1,0].$$
$$\int_{-1}^0|\Lambda^{s+\alpha}v(t,\omega;t_0,\theta_0-z(t_0,\omega))|^2dt\leq c_2(\omega).$$

\proof In the following we will prove some useful estimates in the space of $H^s$ for $s>2(1-\alpha)$ to get an absorbing ball in the space $H^s$.

[\textbf{$L^2$-norm estimates}] First we give the $L^2$-norm estimates which will be used in the proof of the $H^s$-norm estimates. Multiplying (3.6) with $v$ and taking the inner product in $L^2$, we have

$$\aligned\frac{1}{2}\frac{d}{dt}|v|^2+\kappa|\Lambda^\alpha v|^2=&(-u
\cdot\nabla (v+z),v)+(\gamma z,v)\\\leq& C\|\nabla z\|_{L^\infty}[|v|^2+|v|\cdot|z|]+\gamma|z|\cdot|v|.\endaligned$$
Then we obtain
$$\aligned\frac{d}{dt}|v|^2+\frac{\kappa}{2}|\Lambda^\alpha v|^2\leq&[-\lambda_1+c_1\|\nabla z\|_{L^\infty}]|v|^2+c\|\nabla z\|_{L^\infty}^2\cdot|z|^2+c\gamma|z|^2.\endaligned$$
Now we set
$$\mu(t)=-\lambda_1+c_1\|\nabla z(t)\|_{L^\infty},\quad p(t)=c\|\nabla z(t)\|_{L^\infty}^2\cdot|z(t)|^2+c\gamma|z(t)|^2.$$
 Gronwall's lemma yields that
\begin{equation}|v(-2)|^2\leq e^{\int_{t_0}^{-2}\mu(l)dl}|v(t_0)|^2+\int_{t_0}^{-2}e^{\int_\sigma^{-2}\mu(l)dl}p(\sigma)d\sigma.\end{equation}
By (3.2) and (3.4), we can choose $\gamma$ large enough such that $$\lim_{t_0\rightarrow-\infty}\frac{1}{-2-t_0}\int_{t_0}^{-2}(-\lambda_1+c_1\|\nabla z\|_{L^\infty} )dl\leq -\frac{\lambda_1}{4}\quad P-a.s..$$
which combining (3.5) implies that
$$\lim_{t_0\rightarrow-\infty}e^{\int_{t_0}^{-2}\mu(l)dl}=0\quad P-a.s.,$$
and $$\int^{-1}_{-\infty}e^{\int_\sigma^{-2}\mu(l)dl}p(\sigma)d\sigma<\infty \quad P-a.s..$$
By a similar argument as (3.17), we have that for $t\in[-2,-1]$
\begin{equation}|v(t)|^2\leq e^{\int_{-2}^t\mu(l)dl}|v(-2)|^2+\int_{-2}^te^{\int_\sigma^t\mu(l)dl}p(\sigma)d\sigma.\end{equation}
\begin{equation}\int_{-2}^{-1}|\Lambda^\alpha v(l)|^2dl\leq C(|v(-2)|^2+\int_{-2}^{-1}|\mu(l)|dl\sup_{-2\leq t\leq -1}|v(t)|^2+\int_{-2}^{-1}p(l)dl).\end{equation}
Therefore, by (3.17), (3.18) and (3.19) we get that
\begin{equation}\int_{-2}^{-1}|\Lambda^\alpha v(l)|^2dl\leq C(e^{\int_{t_0}^{-2}\mu(l)dl}|v(t_0)|^2\mu_2+p_1),\end{equation}
where
$$\mu_2=1+\int_{-2}^{-1}|\mu(l)|dl\sup_{-2\leq t\leq -1}e^{\int_{-2}^t\mu(l)dl},$$
$$\aligned p_1=&\mu_2\int_{t_0}^{-1}e^{\int_\sigma^{-2}\mu(l)dl}p(\sigma)d\sigma+\int_{-2}^{-1}p(l)dl.\endaligned$$
By (3.2), (3.4), (3.5) the regularity of $z$ and similar arguments as above we have that
$$\sup_{t_0<-1}p_1<\infty \quad P-a.s.$$
[\textbf{$H^s$-norm estimates}] Since $v\in C([t_0,\infty),H^s)\cap L^2_{\rm{loc}}([t_0,\infty),H^{s+\alpha})$, we obtain  the following estimate as (3.14)
\begin{equation}\aligned&\frac{1}{2}\frac{d}{dt}|\Lambda^{\alpha}v|^2+\kappa |\Lambda^{2\alpha} v|^2\\\leq & C|\Lambda^{2\alpha} v||\Lambda (v+z)^2|+\gamma|\Lambda^{\alpha}z||
\Lambda^{\alpha}v|\\\leq & C|\Lambda^{2\alpha} v|[|\Lambda^{1+\sigma_1}v|\|v\|_{L^p}+|\Lambda^{1+\sigma_1}v|\|z\|_{L^p}
+|\Lambda^{1+\tilde{\sigma}_1}z|\|v\|_{L^{\tilde{p}}}\\&+|\Lambda^{1+\sigma_1}z|\|z\|_{L^p}]+\gamma|\Lambda^{\alpha}z||
\Lambda^{\alpha}v|\\\leq & \frac{\kappa}{4} |\Lambda^{2\alpha} v|^2+C(\|v\|_{L^p}^r+\|z\|_{L^p}^r+\varepsilon)|\Lambda^{\alpha} v|^2+|\Lambda^{1+\tilde{\sigma}_1}z||\Lambda^{\alpha}v|^{1-\tilde{r}_0}|\Lambda^{2\alpha}v|^{1+\tilde{r}_0}\\&+C(\gamma|\Lambda^{\alpha}z|^2
+|\Lambda^{1+\sigma_1}z|^2\|z\|_{L^p}^2)\\\leq & \frac{\kappa}{2} |\Lambda^{2\alpha} v|^2+C(\|v\|_{L^p}^r+\|z\|_{L^p}^r+|\Lambda^{1+\tilde{\sigma}_1}z|^{\frac{2}{1-\tilde{r}_0}}+\varepsilon)|\Lambda^{\alpha} v|^2\\&+C(\gamma|\Lambda^{\alpha}z|^2+|\Lambda^{1+\sigma_1}z|^2\|z\|_{L^p}^2),\endaligned\end{equation}
where $r=\frac{2\alpha}{2\alpha-1-\sigma_1}, p=\frac{2}{\sigma_1}$ as in (3.14), $\tilde{r}_0=\frac{1-\tilde{\sigma}_1-\alpha}{\alpha}$ for some $0<\tilde{\sigma}_1=\frac{2}{\tilde{p}}<1-\alpha$ and we used Lemma 2.7 in the second inequality and Lemma 2.8, the interpolation inequality, $H^{2\alpha}\subset L^{\tilde{p}}$ and Young's inequality in the last two inequalities.
Then we get
\begin{equation}\aligned\frac{d}{dt}|\Lambda^{\alpha}v|^2\leq & C(\|v\|_{L^p}^r+\|z\|_{L^p}^r+|\Lambda^{1+\tilde{\sigma}_1}z|^{\frac{2}{1-\tilde{r}_0}})|\Lambda^{\alpha} v|^2\\&+C(\gamma|\Lambda^{\alpha}z|^2+|\Lambda^{1+\sigma_1}z|^2\|z\|_{L^p}^2).\endaligned\end{equation}
By (3.22), (3.16) and Gronwall's lemma, for $l\in[-2,-1]$, we have
\begin{equation}\aligned|\Lambda^\alpha v(-1)|^2\leq & C(|\Lambda^{\alpha}v(l)|^2+\int_{l}^{-1}(\gamma|\Lambda^{\alpha}z|^2+|\Lambda^{1+\sigma_1}z|^2\|z\|_{L^p}^2)d\tau)
\\&\exp\int_{l}^{-1} C[\|v\|_{L^p}^r+\|z\|_{L^p}^r+|\Lambda^{1+\tilde{\sigma}_1}z|^{\frac{2}{1-\tilde{r}_0}}]d\tau \\\leq&C(|\Lambda^{\alpha}v(l)|^2+\int_{-2}^{-1}(\gamma|\Lambda^{\alpha}z|^2
+|\Lambda^{1+\sigma_1}z|^2\|z\|_{L^p}^2)d\tau)\\&\exp\int_{-2}^{-1} C[(\|v(t_0)\|_{L^p}\exp\{-\frac{2\lambda_1}{p}(\tau-t_0)+\int_{t_0}^\tau \|\nabla z(l)\|_{L^\infty} dl\}\\&+\int_{t_0}^\tau(\|\nabla z(l)\|_{L^\infty}\|z(l)\|_{L^p}+C\gamma\|z(l)\|_{L^p})\exp\{-\frac{2\lambda_1}{p}(\tau-l)+\int_{l}^\tau \|\nabla z(\sigma)\|_{L^\infty} d\sigma\}dl)^r\\&+\|z\|_{L^p}^r+|\Lambda^{1+\tilde{\sigma}_1}z|^{\frac{2}{1-\tilde{r}_0}}]d\tau.  \endaligned\end{equation}
Integrating $l$ over $[-2,-1]$ and by (3.20), we obtain
\begin{equation}\aligned|\Lambda^{\alpha}v(-1)|^2\leq&C(\int_{-2}^{-1}|\Lambda^{\alpha}v(l)|^2dl+\int_{-2}^{-1}(\gamma|\Lambda^{\alpha}z|^2
+|\Lambda^{1+\sigma_1}z|^2\|z\|_{L^p}^2)d\tau)\\&\exp\int_{-2}^{-1} C[(\|v(t_0)\|_{L^p}\exp\{-\frac{2\lambda_1}{p}(\tau-t_0)+\int_{t_0}^\tau \|\nabla z(l)\|_{L^\infty} dl\}\\&+\int_{t_0}^\tau(\|\nabla z(l)\|_{L^\infty}\|z(l)\|_{L^p}+C\gamma\|z(l)\|_{L^p})\exp\{-\frac{2\lambda_1}{p}(\tau-l)+\int_{l}^\tau \|\nabla z(\sigma)\|_{L^\infty} d\sigma\}dl)^r\\&+\|z\|_{L^p}^r+|\Lambda^{1+\tilde{\sigma}_1}z|^{\frac{2}{1-\tilde{r}_0}}]d\tau\\\leq& C(e^{\int_{t_0}^{-2}\mu(l)dl}|v(t_0)|^2\mu_2+p_1+p_2)\\&\exp C[\|v(t_0)\|_{L^p}^r\exp\{\frac{r2\lambda_1}{p}t_0+\int_{t_0}^{-1} r\|\nabla z(l)\|_{L^\infty} dl\}+p_3]\\\leq&C(e^{\int_{t_0}^{-2}\mu(l)dl}|v(t_0)|^2\mu_2+p_1+p_2)e^{p_3}\exp [C\|v(t_0)\|^r_{L^p}\exp\{\frac{r2\lambda_1}{p}t_0+\int_{t_0}^{-1} r\|\nabla z(l)\|_{L^\infty} dl\}],\endaligned\end{equation}
where $$p_2=\int_{-3}^{0}(\gamma|\Lambda^{\alpha}z|^2
+|\Lambda^{1+\sigma_1}z|^2\|z\|_{L^p}^2)d\tau,$$
$$\aligned p_3=&C\sup_{t_0<-1}(\int_{t_0}^0(\|\nabla z(l)\|_{L^\infty}\|z(l)\|_{L^p}+C\gamma\|z(l)\|_{L^p})\exp\{\frac{2\lambda_1}{p}l+\int_{l}^0 \|\nabla z(\sigma)\|_{L^\infty} d\sigma\}dl)^r\\&+\int_{-3}^{0}(\|z\|_{L^p}^r+|\Lambda^{1+\tilde{\sigma}_1}z|^{\frac{2}{1-\tilde{r}_0}})d\tau,\endaligned$$

By (3.2), (3.4) (3.5) and similar arguments as above, we can find $\gamma$ large enough and obtain $p_3<\infty$ $P-a.s..$

Moreover, by the same arguments as the proof of (3.23) and (3.24) we have
$$|\Lambda^{\alpha}v(-2)|^2\leq C(e^{\int_{t_0}^{-2}\mu(l)dl}|v(t_0)|^2\mu_1'+p_1'+p_2)e^{p_3}\exp [C\|v(t_0)\|^r_{L^p}\exp\{\frac{r2\lambda_1}{p}t_0+\int_{t_0}^{-1} r\|\nabla z(l)\|_{L^\infty} dl\}],$$
where
$$\mu_1'=1+\int_{-3}^{-1}|\mu(l)|dl\sup_{-3\leq t\leq -1}e^{\int_{-3}^t\mu(l)dl},$$
$$\aligned p_1'=&\mu_1'\int_{t_0}^{-1}e^{\int_\sigma^{-3}\mu(l)dl}p(\sigma)d\sigma+\int_{-3}^{-1}p(l)dl.\endaligned$$
We can easily deduce that $\mu_1'<\infty,\sup_{t_0<-1}p_1'<\infty$ $P$-a.s..
(3.22) yields that for $t\in[-2,-1]$
\begin{equation}\aligned |\Lambda^\alpha v(t)|^2\leq&(|\Lambda^\alpha v(-2)|^2+p_2)e^{p_3}\exp [C\|v(t_0)\|_{L^p}^r\exp\{\frac{r2\lambda_1}{p}t_0+\int_{t_0}^{-1} r\|\nabla z(l)\|_{L^\infty} dl\}]\\\leq&C(e^{\int_{t_0}^{-3}\mu(l)dl}|v(t_0)|^2\mu_1'+p_1'+2p_2)e^{2p_3}\exp [C\|v(t_0)\|_{L^p}^r\exp\{\frac{r2\lambda_1}{p}t_0+\int_{t_0}^{-1} r\|\nabla z(l)\|_{L^\infty} dl\}].\endaligned\end{equation}
Using (3.21) we obtain
\begin{equation}\aligned&\int_{-2}^{-1}|\Lambda^{2\alpha}v(t)|^2dt\\\leq&C[|\Lambda^\alpha v(-2)|^2+\int_{-2}^{-1}(\|v\|_{L^p}^r+\|z\|_{L^p}^r+|\Lambda^{1+\tilde{\sigma}_1}z|^{\frac{2}{1-\tilde{r}_0}})dl\sup_{-2\leq t\leq-1}|\Lambda^\alpha v(t)|^2+p_2]\\\leq&C(e^{\int_{t_0}^{-2}\mu(l)dl}|v(t_0)|^2\mu_1'+3p_2+2p_1')e^{3p_3}\exp [C\|v(t_0)\|_{L^p}^r\exp\{\frac{r2\lambda_1}{p}t_0+\int_{t_0}^{-1} r\|\nabla z(l)\|_{L^\infty} dl\}]+Cp_2,\endaligned\end{equation}
where we used (3.25) in the last inequality.

By a similar argument as (3.14) we have for $s_0>(2-2\alpha)\vee(1-\sigma_0)$
\begin{equation}\aligned&\frac{1}{2}\frac{d}{dt}|\Lambda^{s_0}v|^2+\kappa |\Lambda^{{s_0}+\alpha} v|^2\\\leq & C|\Lambda^{{s_0}+\alpha} v|[|\Lambda^{{s_0}-\alpha+1+\sigma_1}v|\|v\|_{L^p}+|\Lambda^{{s_0}-\alpha+1+\sigma_1}v|\|z\|_{L^p}
+|\Lambda^{{s_0}-\alpha+1+\sigma_1}z|\|v\|_{L^p}\\&+|\Lambda^{{s_0}-\alpha+1+\sigma_1}z|\|z\|_{L^p}]+\gamma|\Lambda^{{s_0}}z||
\Lambda^{s_0}v|\\\leq & \frac{\kappa}{2} |\Lambda^{{s_0}+\alpha} v|^2+C(\|v\|_{L^p}^r+\|z\|_{L^p}^r+|\Lambda^{{s_0}-\alpha+1+\sigma_1}z|^{{2}}+\varepsilon)|\Lambda^{s_0} v|^2\\&+C(\gamma|\Lambda^{{s_0}}z|^2+|\Lambda^{{s_0}-\alpha+1+\sigma_1}z|^2\|z\|_{L^p}^2),\endaligned\end{equation}
where we used Lemmas 2.7, 2.8, the interpolation inequality and Young's inequality in the last two inequalities. Here $r=\frac{2\alpha}{2\alpha-1-\sigma_1}, p=\frac{2}{\sigma_1}$ as in (3.14) and we use $H^{s_0}\subset L^p$.
Therefore  by a similar argument as in the proof of (3.24) and using (3.27) for $s_0=2\alpha$,  we get a similar estimate as (3.24) for $|\Lambda^{2\alpha}v(-1)|^2$.
Thus by  a boot-strapping argument we get that for $s>2(1-\alpha)$
\begin{equation}\aligned|\Lambda^sv(-1)|^2\leq&C(e^{\int_{t_0}^{-2}\mu(l)dl}|v(t_0)|^2\mu_3+q_2)e^{q_3}\exp [C\|v(t_0)\|_{L^p}^r\exp\{\frac{r2\lambda_1}{p}t_0+\int_{t_0}^{-1} r\|\nabla z(l)\|_{L^\infty} dl\}],\endaligned\end{equation}
for suitable $\mu_3, q_2, q_3$. By (3.2), (3.4), (3.5), we can choose $\gamma$ large enough and obtain $\mu_3, q_2, q_3<\infty$ $P-a.s..$
Moreover, we have that
$$\exp\{\frac{r2\lambda_1}{p}t_0+\int_{t_0}^{-1} r\|\nabla z(l)\|_{L^\infty} dl\}\rightarrow0 \textrm{ as } t_0\rightarrow-\infty \qquad P-a.s.,$$
and $$e^{\int_{t_0}^{-2}\mu(l)dl}\rightarrow0 \textrm{ as } t_0\rightarrow-\infty \qquad P-a.s..$$
Then for $|\Lambda^s \theta_0|\leq \rho$, choose $t(\omega)$ such that
$$e^{\int_{t_0}^{-2}\mu(l)dl}|v(t_0)|^2\mu_3\leq 1,$$
\begin{equation}\|v(t_0)\|_{L^p}^r\exp\{\frac{r2\lambda_1}{p}t_0+\int_{t_0}^{-1} r\|\nabla z(l)\|_{L^\infty} dl\}\leq 1,\end{equation}
for all $t_0\leq t(\omega)$, which implies the first result by (3.28).

Furthermore, (3.27) yields that for $t\in[-1,0]$
$$\aligned|\Lambda^s v(t)|^2\leq&C(e^{\int_{t_0}^{-2}\mu(l)dl}|v(t_0)|^2\mu_3+q_4)e^{q_5}\exp [C\|v(t_0)\|_{L^p}^r\exp\{\frac{r2\lambda_1}{p}t_0+\int_{t_0}^{0} r\|\nabla z(l)\|_{L^\infty} dl\}],\endaligned$$
and
$$\aligned&\int_{-1}^{0}|\Lambda^{s+\alpha}v(t)|^2dt\\\leq&C(e^{\int_{t_0}^{-2}\mu(l)dl}|v(t_0)|^2\mu_3+q_6)e^{q_7}\exp [C\|v(t_0)\|_{L^p}^r\exp\{ \frac{r2\lambda_1}{p}t_0+\int_{t_0}^{0} r\|\nabla z(l)\|_{L^\infty} dl\}]+q_8,\endaligned$$
for suitable $q_4,q_5,q_6,q_7,q_8$. By (3.2), (3.4) and (3.5) we can choose $\gamma$ large enough and obtain $q_4,q_5,q_6,q_7,q_8<\infty$ $P$-a.s..

From this and a similar argument as above, the results follow.
$\hfill\Box$
\vskip.10in
\subsection{Compact absorption}

\th{Lemma 3.6} Suppose the condition (E.1) holds. There exists a random radius $r_2(\omega)>0$, such that for all $\rho>0$ there exists $t(\omega)\leq -1$ such that the following holds $P$-a.s. For all $t_0\leq t(\omega)$ and all $\theta_0\in H^s$ with $|\Lambda^s \theta_0|\leq \rho$, the solution $\theta(t,\omega;t_0,\theta_0)=v(t,\omega;t_0,\theta_0-z(t_0,\omega))+z(t,\omega)$ with $v(t_0)=\theta_0-z(t_0,\omega)$ satisfies the following inequality
$$|\Lambda^{s+\delta}\theta(0,\omega;t_0,\theta_0)|^2\leq r_2^2(\omega),$$
for some $0<\delta<\sigma_0\wedge\alpha$.

\proof For $0<\delta<\sigma_0\wedge\alpha$, by Lemma 3.5 we have for almost every $l\in [-1,0]$,  $v(l)\in H^{s+\delta}$. Then by a similar argument as the proof of Theorem 3.1 we obtain the solution $v\in L_{\rm{loc}}^\infty([l,\infty); H^{s+\delta})\cap L^2_{\rm{loc}}([l,\infty);H^{s+\alpha+\delta})$. By a similar estimate as (3.14) we have for  $\sigma_1,r,p$ as in (3.14),
$$\aligned\frac{d}{dt}|\Lambda^{s+\delta}v|^2+\kappa |\Lambda^{s+\alpha+\delta} v|^2\leq & C|\Lambda^{s+\alpha+\delta} v|[|\Lambda^{s+1-\alpha+\delta+\sigma_1}v|\|v\|_{L^p}+|\Lambda^{s+1-\alpha+\delta+\sigma_1}v|\|z\|_{L^p}
\\&+|\Lambda^{s+1-\alpha+\delta+\sigma_1}z|\|v\|_{L^p}+|\Lambda^{s+1-\alpha+\delta+\sigma_1}z|\|z\|_{L^p}]
\\&+C\gamma|\Lambda^{s+\delta}z||\Lambda^{s+\delta}v|\\\leq & \frac{\kappa}{2} |\Lambda^{s+\alpha+\delta} v|^2+C(\|v\|_{L^p}^r+\|z\|_{L^p}^r+|\Lambda^{s+1-\alpha+\delta+\sigma_1}z|^2+\varepsilon)|\Lambda^{s+\delta} v|^2\\&+C(\gamma|\Lambda^{s+\delta}z|^2+|\Lambda^{s+1-\alpha+\delta+\sigma_1}z|^2\|z\|_{L^p}^2),\endaligned$$
where we  choose $\sigma_1$ such that $\sigma_1+\delta<\sigma_0$ and use (E.1) with $z\in C(\mathbb{R};H^{s+1-\alpha+\delta+\sigma_1})$.
Hence by Gronwall's lemma and (3.16) we obtain for $l\in[-1,0]$
$$\aligned|\Lambda^{s+\delta}v(0)|^2\leq&C(|\Lambda^{s+\delta}v(l)|^2+\int_{l}^{0}(\gamma|\Lambda^{s+\delta}z|^2
+|\Lambda^{s+1-\alpha+\delta+\sigma_1}z|^2\|z\|_{L^p}^2)d\tau)\\&\exp\int_{l}^{0} [C(\|v(t_0)\|_{L^p}\exp\{-\frac{2\lambda_1}{p}(\tau-t_0)+\int_{t_0}^\tau \|\nabla z(l_1)\|_{L^\infty} dl_1\}\\&+\int_{t_0}^\tau(\|\nabla z(l_1)\|_{L^\infty}\|z(l_1)\|_{L^p}+C\gamma\|z(l_1)\|_{L^p})\exp\{-\frac{2\lambda_1}{p}(\tau-l_1)+\int_{l_1}^\tau \|\nabla z(\sigma)\|_{L^\infty} d\sigma\}dl_1)^r\\&+\|z\|_{L^p}^r+|\Lambda^{s+1-\alpha+\delta+\sigma_1}z|^2]d\tau.  \endaligned$$
Integrating in $l$ over $[-1,0]$ we deduces that for $\delta\leq\alpha$
$$\aligned|\Lambda^{s+\delta}v(0)|^2\leq&C(\int_{-1}^0|\Lambda^{s+\alpha}v(l)|^2dl+\int_{-1}^{0}(\gamma|\Lambda^{s+\delta}z|^2
+|\Lambda^{s+1-\alpha+\delta+\sigma_1}z|^2\|z\|_{L^p}^2)d\tau)\\&\exp\int_{-1}^{0} [C(\|v(t_0)\|_{L^p}\exp\{-\frac{2\lambda_1}{p}(\tau-t_0)+\int_{t_0}^\tau \|\nabla z(l_1)\|_{L^\infty} dl_1\}\\&+\int_{t_0}^\tau(\|\nabla z(l_1)\|_{L^\infty}\|z(l_1)\|_{L^p}+C\gamma\|z(l_1)\|_{L^p})\exp\{\frac{2\lambda_1}{p}l_1+\int_{l_1}^\tau \|\nabla z(\sigma)\|_{L^\infty} d\sigma\}dl_1)^r\\&+\|z\|_{L^p}^r+|\Lambda^{s+1-\alpha+\delta+\sigma_1}z|^2]d\tau\\\leq&C(\int_{-1}^0
|\Lambda^{s+\alpha}v(l)|^2dl+\int_{-1}^{0}(\gamma|\Lambda^{s+\delta}z|^2
+|\Lambda^{s+1-\alpha+\delta+\sigma_1}z|^2\|z\|_{L^p}^2)d\tau)\\&e^{p_4}\exp [C\|v(t_0)\|^r_{L^p}\exp\{\frac{r2\lambda_1}{p}t_0+\int_{t_0}^{0} r\|\nabla z(l)\|_{L^\infty} dl\}], \endaligned$$
where $$\aligned p_4=&C\sup_{t_0<-1}(\int_{t_0}^0(\|\nabla z(l)\|_{L^\infty}\|z(l)\|_{L^p}+C\gamma\|z(l)\|_{L^p})\exp\{\frac{2\lambda_1}{p}l+\int_{l}^0 \|\nabla z(\sigma)\|_{L^\infty} d\sigma\}dl)^r\\&+\int_{-3}^{0}(\|z\|_{L^p}^r+|\Lambda^{s+1-\alpha+\delta+\sigma_1}z|^2)d\tau,\endaligned$$
 By (3.2), (3.4), (3.5) we know $p_4<\infty$ $P$-a.s. which combining Lemma 3.5 and (3.29) implies the absorption of $\varphi$ in $H^{s+\delta}$ at time $t =0$.
$\hfill\Box$

\vskip.10in

Since the embedding $ H^{s+\delta}\subset H^s$ is compact, by Proposition 2.4 and [8, Corollary 4.6] we obtain the following results.
\vskip.10in
\th{Theorem 3.7} Fix $\alpha>1/2$. Suppose the condition (E.1) holds. Then the stochastic flow associated with the quasi-geostrophic equation (3.1) driven by additive  noise has a compact stochastic attractor in $H^s$.

Moreover, the Markov semigroup induced by the flow on $H^s$ has an invariant measure $\rho$.

\section{Upper semicontinuity of random attractors}
In this section we consider the following equation
\begin{equation}d\theta+(A_\alpha\theta+u\cdot\nabla\theta)dt=\varepsilon dW.\end{equation}
 Now we fix the same $s$ as in Section 3 and assume that $G$ satisfies  (E.1). By [13, Theorem 5.1], the solution operator $S:S(t)\theta_0=\theta(t,\theta_0)$ defines a semigroup in the space $H^s$, where $\theta(t,\theta_0)$ is the solution of equation (4.1) with $\varepsilon=0$ and initial value $\theta_0$ at time $0$. Moreover, $\{S(t)\}$ possesses a global attractor $\mathcal{A}$ in $H^s$.

By Theorem 3.4 we obtain a continuous  random dynamical system associated with (4.1)
$$\varphi_\varepsilon:\mathbb{R}^+\times \Omega\times H^s\rightarrow H^s.$$

First we prove for $P$-a.e. $\omega\in \Omega$ and $\theta_0\in H^s, t_0\in \mathbb{R}^+$
$$\varphi_\varepsilon(t_0,\vartheta_{-t_0}\omega)\theta_0\rightarrow S(t_0)\theta_0 \textrm{ as } \varepsilon\rightarrow0,$$
i.e.
$$S_\varepsilon(0,-t_0;\omega)\theta_0\rightarrow S(t_0)\theta_0 \textrm{ as } \varepsilon\rightarrow0,$$
where $S_\varepsilon(0,-t_0;\omega)\theta_0$ denote the stochastic flow associated with equation (4.1) obtained in Section 3.

\th{Proposition 4.1} Suppose the condition (E.1) holds.  Then for $P$-a.e. $\omega\in \Omega$ and $t_0\in \mathbb{R}^+$ and $B\subset H^s$ bounded
$$\lim_{\varepsilon\rightarrow0}\sup_{\theta_0\in B}|\Lambda^s[S_\varepsilon(0,-t_0,\omega)\theta_0-\theta(t_0;\theta_0)]|=0.$$
\proof Denote $\theta_\varepsilon(t,\omega)=S_\varepsilon(t,-t_0;\omega)\theta_0$ for simplicity. Let $\zeta_\varepsilon(t,\omega)=\theta_\varepsilon(t,\omega)-\theta(t)$ where $\theta(t)$ is the solution to the unperturbed equations with the same initial condition $\theta_0$ at $-t_0$. Then $\zeta_\varepsilon$ satisfies
$$d\zeta_\varepsilon+A_\alpha\zeta_\varepsilon dt+(u_{\zeta_\varepsilon}\cdot\nabla\zeta_\varepsilon+u_{\zeta_\varepsilon}\cdot\nabla\theta+u_\theta\cdot\nabla\zeta_\varepsilon)dt=\varepsilon   dW(t),$$where $u_{\zeta_\varepsilon}$ satisfies (1.3) with $\theta$ replaced by ${\zeta_\varepsilon}$.
We use the change of variable
$$\eta_\varepsilon=\zeta_\varepsilon-z_\varepsilon:=\zeta_\varepsilon-\varepsilon\int_{-t_0}^te^{(t-l)A_\alpha}dW(l),$$
 which satisfies the following equality in the weak sense,
$$\frac{d\eta_\varepsilon}{dt}+A_\alpha\eta_\varepsilon+u_{\eta_\varepsilon+z_\varepsilon}\cdot\nabla(\eta_\varepsilon+z_\varepsilon)+u_{\eta_\varepsilon+z_\varepsilon}
\cdot\nabla\theta+u_\theta\cdot \nabla(\eta_\varepsilon+z_\varepsilon)=0,$$
where $u_{\eta_\varepsilon+z_\varepsilon}=u_{\zeta_\varepsilon}$.
Since $\theta,\theta_\varepsilon\in C([-t_0,T],H^s)\cap L^2([-t_0,T],H^{s+\alpha})$,  we obtain the following estimate  as in (3.14) by taking the scalar product with $\Lambda^{s}e_k$, multiplying both sides by $\langle \eta_\varepsilon,\Lambda^se_k\rangle$ and summing up over $k$, ,
$$\aligned&\frac{1}{2}\frac{d}{dt}|\Lambda^s\eta_\varepsilon|^2+\kappa|\Lambda^{s+\alpha}\eta_\varepsilon|^2
\\\leq& \langle u_{\eta_\varepsilon+z_\varepsilon}\cdot\nabla(\eta_\varepsilon+z_\varepsilon)+u_{\eta_\varepsilon+z_\varepsilon}\cdot\nabla\theta+u_\theta\cdot \nabla(\eta_\varepsilon+z_\varepsilon),\Lambda^{2s}\eta_\varepsilon\rangle\\\leq&
C(\|\eta_\varepsilon\|_{L^p}^r|\Lambda^s\eta_\varepsilon|^2
+\|z_\varepsilon\|_{L^p}^r|\Lambda^{s}\eta_\varepsilon|^2+|\Lambda^{s-\alpha+1+\sigma_1}z_\varepsilon|^2
|\Lambda^{s}\eta_\varepsilon|^2+\|z_\varepsilon\|_{L^p}^2|\Lambda^{s-\alpha+1+\sigma_1}z_\varepsilon|^2)\\&+C[\|\theta\|_{L^p}^r+|\Lambda^{s-\alpha+1+\sigma_1}\theta|^2]|
\Lambda^{s}\eta_\varepsilon|^2+C|\Lambda^{s-\alpha+1+\sigma_1}z_\varepsilon|^2\|\theta\|_{L^p}^2+C|\Lambda^{s+\alpha}\theta|^2\|
z_\varepsilon\|_{L^p}^2\\&+\frac{\kappa}{2}|\Lambda^{s+\alpha}\eta_\varepsilon|^2,\endaligned$$
where $\sigma_1, r, p$ are as in (3.14)  and we used Lemmas 2.7, 2.8, the interpolation inequality and Young's inequality in the last inequality. Here the calculation in the last inequality is similar as in (3.14), so we omit the details. Then we have
$$\aligned&\frac{1}{2}\frac{d}{dt}|\Lambda^s\eta_\varepsilon|^2\leq h(t)+k(t)|
\Lambda^{s}\eta_\varepsilon|^2,\endaligned$$
where $$h(t)= C(|\Lambda^{s-\alpha+1+\sigma_1}z_\varepsilon|^2\|\theta\|_{L^p}^2+|\Lambda^{s+\alpha}\theta|^2\|z_\varepsilon\|_{L^p}^2
+\|z_\varepsilon\|_{L^p}^2|\Lambda^{s-\alpha+1+\sigma_1}z_\varepsilon|^2).$$
$$k(t)=C[\|\eta_\varepsilon\|_{L^p}^r+\|z_\varepsilon\|_{L^p}^r+\|\theta\|_{L^p}^r+|\Lambda^{s-\alpha+1+\sigma_1}\theta|^2+|\Lambda^{s-\alpha+1+\sigma_1}z_\varepsilon|^2].$$
Since $\theta\in C([-t_0,\infty),H^s)\cap L^2_{\rm{loc}}([-t_0,\infty),H^{s+\alpha})$ and $\sup_{\theta_0\in B}h(t)\rightarrow0$ when $\varepsilon\rightarrow0$, by Gronwall's lemma we obtain
$$\sup_{\theta_0\in B}|\Lambda^s \eta_\varepsilon(t)|^2\rightarrow0 \textrm{ as }\varepsilon\rightarrow0,$$
for all $t\geq -t_0$.
Therefore $$\sup_{\theta_0\in B}|\Lambda^s \zeta_\varepsilon(t)|^2\rightarrow0 \textrm{ as }\varepsilon\rightarrow0.$$
$\hfill\Box$
\vskip.10in
By the computation in Lemma 3.6 we can easily check that
$$\lim_{\varepsilon\rightarrow0}r_{2,\varepsilon}(\omega)\leq r_d,$$
with $r_d$ independent of $\omega\in \Omega$, where $r_{2,\varepsilon}$ is the random radius for the solution to (4.1) we obtained in Lemma 3.6.

Then by [5, Theorem 2, Lemma 1] we obtain
\vskip.10in
\th{Theorem 4.2} Suppose the condition (E.1) holds. Let $\mathcal{A}_\varepsilon(\omega)$ denote the random attractor  for $\varphi_\varepsilon$. Then
$$\lim_{\varepsilon\rightarrow0}d(\mathcal{A}_\varepsilon,\mathcal{A})=0 \qquad P-a.s.$$
Moreover, the convergence above is upper semicontinuous in $\varepsilon$, that is
$$\lim_{\varepsilon\rightarrow\varepsilon_0}d(\mathcal{A}_\varepsilon(\omega),\mathcal{A}_{\varepsilon_0}(\omega))=0 \qquad P-a.s.$$

\section{The triviality of the random attractor}

In this section we assume that $G$ satisfies the same condition as in Section 3 and we take $s=1$ for simplicity. Then our assumption for $G$ is
$$\mathcal{E}_0:=\rm{Tr}(\Lambda^{ (4+2\varepsilon_0)} GG^*)<\infty.$$

Under this condition we will prove that if the viscosity constant $\kappa$ is large enough or $\mathcal{E}_0$ is small enough, the random attractor is trivial. The idea for the proof is inspired by the approach in \cite{M99}. But for the stochastic quasi-geostrophic equation we need more delicate estimates.  Since for the stochastic quasi-geostrophic equation the dissipation term is not strong enough and cannot control the nonlinear term,
 we will use $L^p$-norm estimate to control the nonlinear term in a larger space.
  In the following we will prove for almost every realization of the noise, trajectories starting from different initial conditions in $H^1$ converge to each other in a larger space $H^{-1/2}$ which is the dual space of Sobolev space $H^{1/2}$.
\vskip.10in

\th{Lemma 5.1}Fix $\alpha>1/2$. Suppose the condition (E.1) holds with $s=1$. If $\delta_0=\kappa-2^{p/2}C_R^p C_S^{2p}\kappa^{1-p}[p(p-1)] ^{p/2}\lambda_1^{-p/2}\mathcal{E}_0^{{p} /2}>0$, i.e. $\kappa^{\frac{3}{2}p}>2^{p/2}C_R^p C_S^{2p}[p(p-1)] ^{p/2}\mathcal{E}_0^{{p} /2}$ for $p=\frac{\alpha+1}{\alpha-\frac{1}{2}}$, where $C_S, C_R$ are the constants for Sobolev embedding and Riesz transform respectively, then for $\delta\in(0,\delta_0)$ and $\theta_0\in H^1$, there exists a positive random time $\tau=\tau(t_0,\omega,\theta_0)$ independent of $\tilde{\theta}_0$ such that for all $t>\tau+t_0$
$$|\Lambda^{-1/2}(S(t,t_0;\omega)\theta_0-S(t,t_0;\omega)\tilde{\theta}_0)|^2\leq |\Lambda^{-1/2}(\theta_0-\tilde{\theta}_0)|^2e^{-\delta (t-t_0)}.$$
Moreover, $E\tau^q<\infty$ for any $q\in (0,+\infty)$.

\proof  We obtain that $\rho:=S(\cdot,t_0;\omega)\tilde{\theta}_0-S(\cdot,t_0;\omega){\theta}_0$ satisfies the following equation in the weak sense:
$$\aligned\frac{d\rho(t)}{dt}=&-A_\alpha\rho-\tilde{u}\cdot\nabla\tilde{\theta}+u\cdot \nabla\theta
\\=&-A_\alpha\rho-u\cdot\nabla\rho-u_\rho\cdot\nabla\tilde{\theta},\endaligned$$
where $u_\rho,\tilde{u}$ satisfy (1.3) with $\theta$ replaced by $\rho, \tilde{\theta}$ respectively and we write $\theta=S(\cdot,t_0;\omega)\theta_0, \tilde{\theta}=S(\cdot,t_0;\omega)\tilde{\theta}_0$ for simplicity.
Taking the inner product with $\Lambda^{-1}\rho$ in $H$, and by  $${ }_{H^{-1}}\!\langle u_\rho\cdot\nabla\tilde{\theta}, \Lambda^{-1}\rho\rangle_{H^1}=0,$$
(cf.\cite{Re95}),
we have
$$\frac{1}{2}\frac{d}{dt}|\Lambda^{-\frac{1}{2}}\rho|^2=-\kappa|\Lambda^{\alpha-\frac{1}{2}} \rho|^2-{ }_{H^{-1}}\!\langle u\cdot\nabla\rho, \Lambda^{-1}\rho\rangle_{H^1}.$$

We calculate
$$\aligned|{ }_{H^{-1}}\!\langle u\cdot\nabla\rho, \Lambda^{-1}\rho\rangle_{H^1}|\leq & \|u\|_{L^{p} }\|\rho\|_{L^{p_1}  }\|\nabla \Lambda^{-1}\rho\|_{L^{p_1} }\leq C_S\|u\|_{L^{p} }\|\rho\|_{H^{1/{p} }}\|\nabla \Lambda^{-1}\rho\|_{H^{1/{p} }}\\\leq &C_SC_R\|\theta\|_{L^{p} }\|\Lambda^{-1}\rho\|_{H^{1+\frac{1}{{p} }}}^2\leq C_SC_R\|\theta\|_{L^{p} }\|\Lambda^{-1}\rho\|^{2/r}_{H^{\frac{1}{2}}}\|\Lambda^{-1}\rho\|^{2(1-\frac{1}{r})}_{H^{\frac{1}{2}+\alpha}}\\\leq& \frac{\kappa}{2}|\Lambda^{\alpha-\frac{1}{2}}\rho|^2+ C^r(\frac{\kappa}{2})^{1-r}\|\theta\|_{L^{p} }^r|\Lambda^{-\frac{1}{2}}\rho|^2,\endaligned$$
where $C_S, C_R$ are the constants for Sobolev embedding and Riesz transform, respectively and $C=C_SC_R$ and we used $H^{1/{p} }\subset L^{p_1} $  in the second inequality  and the interpolation inequality in the forth inequality and Young's inequality in the last inequality. Here $\frac{1}{{p} }+\frac{2}{{p_1} }=1$ for $p>\frac{1}{\alpha-\frac{1}{2}}, r=\frac{\alpha}{\alpha-\frac{1}{2}-\frac{1}{{p} }} $. Then we obtain
$$\frac{d}{dt}|\Lambda^{-\frac{1}{2}}\rho|^2\leq-{\kappa}|\Lambda^{\alpha-\frac{1}{2}} \rho|^2+2C^r(\frac{\kappa}{2})^{1-r}\|\theta\|_{L^{p} }^r|\Lambda^{-\frac{1}{2}}\rho|^2.$$
Thus Gronwall's lemma yields that
$$|\Lambda^{-\frac{1}{2}}\rho(t)|^2\leq e^{(t-t_0)\Gamma(t-t_0;t_0,\theta_0)}|\Lambda^{-\frac{1}{2}}\rho(t_0)|^2,$$
where
$$\Gamma(t_1;t_0,\theta_0)=-\kappa+2C^r(\frac{\kappa}{2})^{1-r}\frac{1}{t_1}\int_{t_0}^{t_1+t_0}\|\theta(s)\|_{L^{p} }^rds.$$

By  Proposition A.2 we obtain
$$\aligned&\|\theta(t_0+t_1)\|_{L^{p} }^{p} +\lambda_1\int_{t_0}^{t_0+t_1}\int_{\mathbb{T}^2}|\theta(l)|^{p} d\xi dl\\\leq&\|\theta_0\|_{L^{p} }^{p}
+C_S^p[\frac{1}{2}p(p-1)]^{p/2}\lambda_1^{-\frac{p-2}{2}}\mathcal{E}_0^{p/2}t_1
+p\int_{t_0}^{t_0+t_1}
\int_{\mathbb{T}^2}|\theta(l)|^{p-2}\theta(l) d\xi dW(l).\endaligned$$

Since $p=\frac{\alpha+1}{\alpha-\frac{1}{2}}$ implies $p=r$, we obtain that $$\aligned\Gamma(t_1,t_0;\theta_0) \leq&-\kappa+2C^p(\frac{\kappa}{2})^{1-p}\frac{1}{t_1}\int_{t_0}^{t_0+t_1}\|\theta(s)\|_{L^{p} }^{p} ds\\\leq &-\kappa+2C^p(\frac{\kappa}{2})^{1-p}\frac{1}{t_1\lambda_1}\|\theta_0\|_{L^{p} }^{p}
+2^{p/2}C^p C_S^p\kappa^{1-p}[p(p-1)] ^{p/2}\lambda_1^{-p/2}\mathcal{E}_0^{{p} /2}
\\&+2C^p(\frac{\kappa}{2})^{1-p}\frac{p}{t_1\lambda_1}\int_{t_0}^{t_0+t_1}
\int_{\mathbb{T}^2}|\theta(l)|^{{p} -2}\theta(l) d\xi dW(l).\endaligned$$
For $M(t_1;t_0):={p} \int_{t_0}^{t_0+t_1}
\int_{\mathbb{T}^2}|\theta(l)|^{{p} -2}\theta(l) d\xi dW(l)$, we have
$$\langle M\rangle_{t_1}\leq C{p} ^2\mathcal{E}_0\int_{t_0}^{t_0+t_1}(\int_{\mathbb{T}^2}|\theta(s)|^{{p} -1}d\xi)^2ds,$$
where we use $\|[\sum_j(G(e_j))^2]^{1/2}\|_{L^\infty}\leq [\sum_j\|G(e_j)\|_{L^\infty}^2]^{1/2}\leq C(\sum_j|\Lambda^{1+\varepsilon}G(e_j)|^2)^{1/2}.$
Then for any $m>2$
$$\langle M\rangle_{t_1}^m\leq C{p} ^{2m}\mathcal{E}_0^m (\int_{t_0}^{t_0+t_1}(\int_{\mathbb{T}^2}|\theta(s)|^{{p} -1}d\xi)^2ds)^m\leq C{p} ^{2m}\mathcal{E}_0^mt_1^{m-1}\int_{t_0}^{t_0+t_1}(\int_{\mathbb{T}^2}|\theta(s)|^{2m(p -1)}d\xi)ds.$$
Since $C_0:=\|\theta_0\|_{L^{2m({p} -1)}}^{2m(p-1)}\leq C\|\theta_0\|_{H^1}^{2m(p-1)}<\infty$, by  Proposition A.2 there exists a constant $C_{{p} ,m}(C_0)$ such that $E\|\theta(t)\|_{L^{2m({p} -1)}}^{2m({p} -1)}\leq  C_{{p} ,m}$ for $t\geq t_0$. Thus for $M_n=\sup_{n-1\leq t<n}M(t;t_0)$, we have
$$P(|M_n|>\frac{\varepsilon\lambda_1}{4C^p(\frac{\kappa}{2})^{1-p}} n)\leq \frac{{p} ^{2m}\mathcal{E}_0^mC_{{p} ,m}n^m}{(\frac{\varepsilon\kappa^{p-1}\lambda_1}{2^{p+1}C^p})^{2m}n^{2m}}.$$

Now define the following random times
$$T_{\rm{bound}}(t_0,\omega,\theta_0):=\sup\{n: |M_n|>\frac{\varepsilon\lambda_1}{4C^p(\frac{\kappa}{2})^{1-p}} n\}.$$

 By   Lemma A.3, we have if $m>1$, then $T_{\rm{bound}}$ is finite almost surely.
Define
$$N_1:=\frac{2^{p+1}C^p\|\theta_0\|_{L^p}^p}{\kappa^{p-1}\lambda_1\varepsilon}.$$
 Set
$$\tau=\max(T_{\rm{bound}}, N_1),$$ then we get that
$$n>\tau\Rightarrow \Gamma(n;t_0,\theta_0)-(-\delta_0)<\varepsilon.$$ Now  we obtain for $\delta\in(0,\delta_0)$ and $t>\tau+t_0$,
$$|\Lambda^{-1/2}\rho(t)|^2\leq |\Lambda^{-1/2}(\theta_0-\tilde{\theta}_0)|^2e^{-\delta (t-t_0)}.$$
For $p_0\in(0,+\infty)$ by  Lemma A.3 $E\tau^{p_0}$ is finite.
$\hfill\Box$
\vskip.10in
Now we will prove the main result of this section. First we will prove the existence of the limit of the stochastic flow $S(t,r,\omega)\theta_0$ constructed in Section 3 when time $r$ goes to $-\infty$. Then selecting a strictly stationarity version of the limiting process is the random attractor desired.
\vskip.10in
\th{Theorem 5.2}Fix $\alpha>1/2$. Suppose the condition (E.1) holds with $s=1$. If $\delta_0=\kappa-2^{p/2}C_R^p C_S^{2p}\kappa^{1-p}[p(p-1)] ^{p/2}\lambda_1^{-p/2}\mathcal{E}_0^{{p} /2}>0$ for $p=\frac{\alpha+1}{\alpha-\frac{1}{2}}$, where $C_S, C_R$ are the constants for Sobolev embedding and Riesz transform respectively, then the RDS $\varphi$ associated with the stochastic quasi-geostrophic equation (3.1) has a compact random
attractor $\mathcal{A}(\omega)$ consisting of a single point:
$$\mathcal{A}(\omega)=\{\tilde{\eta}_0(\omega)\}.$$
Moreover,  the invariant measure is unique.

\proof First we prove that for all $t\in\mathbb{R}$ and there exists $\Omega_1\subset \Omega$ such that $P(\Omega_1)=1$ and for $\omega\in\Omega_1$ there exists a limit $\eta_t(\omega)$ such that
$$\lim_{r\rightarrow-\infty}|\Lambda^{-1/2}(S(t,r;\omega)\theta_0-\eta_t(\omega))|=0.$$
For fixed $t_1\in \mathbb{Z}$ define $$n^*(\omega):=\sup\{n:\tau(-n+t_1,\omega,\theta_0)>n\},$$
where $\tau(-n+t_1,\omega,\theta_0)$ is the random time we obtained in Lemma 5.1. By the estimate for $M_n$ in the proof of Lemma 5.1 and definition of $\tau(-n+t_1,\omega,\theta_0)$ we have for $p_0\in (1,+\infty)$
$$E\tau(-n+t_1,\omega,\theta_0)^{p_0}\leq C(p_0)$$
Then by  Lemma A.3 we know
that $E(n^*)^{p_1}<\infty$ for any $p_1\in(1,+\infty)$ and for $t>t_1, n_1,n_2\in \mathbb{Z}^+, n_1>n_2>n^*$ and $\theta_0\in H^1$,
\begin{equation}\aligned&|\Lambda^{-1/2}(S(t,-n_1+t_1;\omega)\theta_0-S(t,-n_2+t_1;\omega)\theta_0))|
\\\leq&\sum_{j=-n_1}^{-n_2-1}|\Lambda^{-1/2}(S(t,j+1+t_1;\omega)\theta_0-S(t,j+t_1;\omega)\theta_0)|
\\=&\sum_{j=-n_1}^{-n_2-1}|\Lambda^{-1/2}(S(t,j+1+t_1;\omega)S(j+t_1+1,j+t_1;\omega)\theta_0-S(t,j+1+t_1;\omega)\theta_0)|\\\leq& \sum_{j=-n_1}^{-n_2-1}|\Lambda^{-1/2}(S(j+1+t_1,j+t_1;\omega)\theta_0-\theta_0)|e^{-\frac{\delta}{2}(t-(j+1+t_1))}
,\endaligned\end{equation}
where we used Lemma 5.1 and $t-(j+1+t_1)>\tau(j+1+t_1,\omega,\theta_0)$ for any $j\in\mathbb{N}, -n_1\leq j\leq -n_2-1$ in the last inequality. Now define the following random time
$$\tau_0:=\sup\{n:|\Lambda^{-1/2}S(-n+1+t_1,-n+t_1;\omega)\theta_0|>\frac{\varepsilon\delta^2}{8}n\},$$
Since $\theta_0\in H^1$ by  Proposition A.2 there exists a constant $C(m)$ such that $$E|S(-n+1+t_1,-n+t_1;\omega)\theta_0|^m\leq  C(m),$$ for $m,n\in\mathbb{N}$. Then by  Lemma A.3 $E\tau_0^p<\infty$ for any $p\in (1,+\infty)$. Define $\Omega_0:=\{\omega:\tau_0(\omega)\vee n^*(\omega)<\infty\}$ and then $P(\Omega_0)=1$.
For $\omega\in \Omega_0$ and $n_1>n_2>n^*\vee\tau_0$, by (5.1) we have
$$|\Lambda^{-1/2}(S(t,-n_1+t_1;\omega)\theta_0-S(t,-n_2+t_1;\omega)\theta_0)|
\leq C\varepsilon e^{-\frac{\delta}{2} t}
.$$

 Therefore,  for all $t>t_1$  there exists a process $\eta_t(\omega)$ such that
\begin{equation}\lim_{n\rightarrow-\infty}|\Lambda^{-1/2}(S(t,n+t_1;\omega)\theta_0-\eta_t(\omega))|=0.\end{equation}
Since $t_1$ is arbitrary we can define $\eta$ for all time. Now we want to prove the convergence in (5.2) is satisfied from any initial time.

For $\omega\in \Omega_0$, $n>n^*\vee\tau_0$ and $r\in [-n-1,-n]$ we obtain for $\tilde{\theta}_0\in H^1$
\begin{equation}\aligned &|\Lambda^{-1/2}(S(t,-n+t_1;\omega)\theta_0-S(t,r+t_1;\omega)\tilde{\theta}_0)|^2\\=
&|\Lambda^{-1/2}(S(t,-n+t_1;\omega)S(-n+t_1,r+t_1;\omega)\tilde{\theta}_0-S(t,-n+t_1;\omega)\theta_0)|^2
\\\leq&(|S(-n+t_1,r+t_1;\omega)\tilde{\theta}_0|^2+C)e^{-\delta(t+n-t_1)}, \endaligned\end{equation}
where we used Lemma 5.1 in the last inequality. By the same argument as (3.17) we have for $r\in[-n-1,-n]$
$$\aligned &|S(-n+t_1,r+t_1;\omega)\tilde{\theta}_0|^2\leq 2|v(-n+t_1,r+t_1,\omega,\tilde{\theta}_0-z(r+t_1))|^2+2|z(-n+t_1)|^2
\\\leq& 2e^{\int_{r+t_1}^{-n+t_1}\mu(l)dl}|\tilde{\theta}_0-z(r+t_1)|^2+2\int_{r+t_1}^{-n+t_1}e^{\int_\sigma^{-n+t_1}\mu(l)dl}p(\sigma)d\sigma+2|z(-n+t_1)|^2,\endaligned$$
where
$$\mu(t)=-\lambda_1+c_1\|\nabla z(t)\|_{L^\infty},\quad p(t)=c\|\nabla z(t)\|_{L^\infty}^2\cdot|z(t)|^2+c\gamma|z(t)|^2.$$
By (3.5)  there exists $\Omega_2\subset \Omega$ such that $P(\Omega_2)=1$ and for any $\varepsilon>0, \omega\in \Omega_2$ there exists $N_0(\omega)$ such that for $t>N_0$, $c_1\|\nabla z(-t)\|_{L^\infty}\leq \varepsilon t.$ By this we obtain for $n>N_0+t_1$
\begin{equation}\aligned |S(-n+t_1,r+t_1;\omega)\tilde{\theta}_0|^2\leq& 2e^{\int_{r+t_1}^{-n+t_1}(-\lambda_1-\varepsilon l)dl}|\tilde{\theta}_0-z(r+t_1)|^2\\&+2\int_{r+t_1}^{-n+t_1}e^{\int_\sigma^{-n+t_1}(-\lambda_1-\varepsilon l)dl}p(\sigma)d\sigma+2|z(-n+t_1)|^2\\\leq& 2e^{\varepsilon (n+\frac{1}{2}-t_1) }|\tilde{\theta}_0-z(r+t_1)|^2\\&+2\int_{-n-1+t_1}^{-n+t_1}e^{\varepsilon (n+\frac{1}{2}-t_1)}p(\sigma)d\sigma+2|z(-n+t_1)|^2
\endaligned\end{equation}
Combining (5.3) and (5.4) we obtain for $n>(n^*\vee\tau_0)\vee(N_0+t_1)$
$$\aligned &\sup_{r\in[-n-1,-n]}|\Lambda^{-1/2}(S(t,-n+t_1;\omega)\theta_0-S(t,r+t_1;\omega)\tilde{\theta}_0)|^2\\\leq& [2e^{\varepsilon (n+\frac{1}{2}-t_1) }|\tilde{\theta}_0-z(r+t_1)|^2\\&+2\int_{-n-1+t_1}^{-n+t_1}e^{\varepsilon (n+\frac{1}{2}-t_1)}p(\sigma)d\sigma+2|z(-n+t_1)|^2+C]e^{-2\delta(t+n-1-t_1)}\endaligned$$
Choosing $\varepsilon<\delta$ we obtain
that  for $\omega\in\Omega_1:=\Omega_0\cap\Omega_2$,
$$\lim_{r\rightarrow-\infty}|\Lambda^{-1/2}(S(t,r+t_1;\omega)\theta_0-\eta_t(\omega))|=0.$$
Moreover, for $\omega\in\Omega_1$, $-\infty<t_1<t_2<+\infty$ and any bounded set $B$ in $H^1$
\begin{equation}\lim_{r\rightarrow-\infty}\sup_{t\in[t_1,t_2]}\sup_{\theta_0\in B}|\Lambda^{-1/2}(S(t,r;\omega)\theta_0-\eta_t(\omega))|=0,\end{equation}
which also implies that $\eta_t(\omega)$ is independent of $\theta_0, t_1$.

For $\theta_0\in H^1$, by similar arguments as  the proof of Lemma 3.6 we obtain that there exists $K(t_1,t_2,\omega)$ such that
 \begin{equation}\sup_{t\in[t_1,t_2]}|\Lambda^{1+\delta}\eta_t(\omega)|\leq \limsup_{r\rightarrow-\infty}\sup_{t\in[t_1,t_2]}|\Lambda^{1+\delta}S(t,r;\omega)\theta_0|\leq K(t_1,t_2,\omega)\quad P-a.s.\end{equation}
Thus  the interpolation inequality and (5.5), (5.6) yield that
\begin{equation}\aligned&\lim_{r\rightarrow-\infty}\sup_{t\in[t_1,t_2]}|\Lambda(S(t,r;\omega)\theta_0-\eta_t(\omega))|\\\leq&\lim_{r\rightarrow-\infty}C
\sup_{t\in[t_1,t_2]}|\Lambda^{-{1/2}}(S(t,r;\omega) \theta_0-\eta_t(\omega))|^{\beta_1}\sup_{t\in[t_1,t_2]}|\Lambda^{1+\delta}(S(t,r;\omega)\theta_0-\eta_t(\omega))|^{1-\beta_1}=0,\endaligned\end{equation}
where $\beta_1=\frac{\delta}{\frac{3}{2}+\delta}$. From (5.6) and (5.7) we know that for almost all $\omega$, $\eta(\omega)\in C(\mathbb{R};H^1)\cap L^2_{\rm{loc}}(\mathbb{R};H^{1+\delta}).$
Since for $t\geq r$
$$S(0,r-t;\vartheta_t\omega)\theta_0=S(t,r;\omega)\theta_0,$$
letting $r\rightarrow-\infty$ and  by (5.5) we obtain
$$\eta_0(\vartheta_t\omega)=\eta_t(\omega)\quad P-a.s.,$$
with the zero set  depending on $t$. Now   we can use  Proposition 2.6 to deduce the existence of an indistinguishable process $\tilde{\eta}$ such that for all $\omega\in\Omega$ $$\tilde{\eta}_0(\vartheta_t\omega)=\tilde{\eta}_t(\omega),$$ and $\tilde{\eta}(\omega)\in C(\mathbb{R};H^1)\cap L^2_{\rm{loc}}(\mathbb{R};H^{1+\delta}).$

Now we define
$$\mathcal{A}(\omega)=\{\tilde{\eta}_0(\omega)\}.$$
Since $x\mapsto\varphi(t,\omega)x$ is continuous in $H^1$, we get that $P$-a.s.
$$\aligned\varphi(t,\omega)\mathcal{A}(\omega)
=&\varphi(t,\omega)\lim_{r\rightarrow-\infty}S(0,r;\omega)\theta_0
\\=&\lim_{r\rightarrow-\infty}S(t, r;\omega)\theta_0\\=&\{\tilde{\eta}_t(\omega)\}\\=&\mathcal{A}(\vartheta_t(\omega)),\endaligned$$
which implies the invariance of $\mathcal{A}$.
Now for any bounded set $B\subset H^1$ $P$-a.s.
$$\aligned&\lim_{r\rightarrow-\infty}\sup_{\theta_0\in B}|\Lambda(S(0,r;\omega)\theta_0-\tilde{\eta}_0(\omega))|\\\leq& C\lim_{r\rightarrow-\infty}\sup_{\theta_0\in B}|\Lambda^{-1/2}(S(0,r;\omega)\theta_0-\tilde{\eta}_0(\omega))|^{\beta_1}|\Lambda^{\delta+1}(S(0,r;\omega) \theta_0-\tilde{\eta}_0(\omega))|^{1-\beta_1}=0,\endaligned$$
which implies that $\mathcal{A}$ attracts all the deterministic bounded sets. Now the first result follows. The uniqueness of the invariant measures is then obvious. $\hfill\Box$
\vskip.10in

\section{Multiplicative noise}

In this section we consider the abstract stochastic evolution equation with Stratonovich multiplicative noise in place of Eqs (1.1)-(1.3),
\begin{equation}d\theta +A_\alpha\theta+u(t)\cdot \nabla\theta(t)dt+\sum_{j=1}^mb_j\theta\circ dw_j(t)=0,\end{equation}
where  $u$ satisfies (1.3), $b_1,...,b_m\in\mathbb{R}$ and  $W=(w_j(t), 1\leq j\leq m)$, are  two-sided  Wiener processes on the canonical Wiener space $(\Omega,\mathcal{F},\mathcal{F}_t,P)$, i.e. $\Omega= C_0(\mathbb{R},\mathbb{R}^m):=\{w\in C(\mathbb{R},\mathbb{R}^m), w(0)=0\}$, $W(\omega)(t):=\omega(t)$, $\mathcal{F}_t$ is canonical filtration and $\vartheta_t$ is the Wiener shift given by $\vartheta_t\omega:=\omega(t+\cdot)-\omega(t)$ and $P=$ the law of $W$.  Here we have that $w_j, 1\leq j\leq m$  have strictly stationary increments, i.e. for all $t,r\in \mathbb{R}$, $\omega\in \Omega$,
$$w_j(t,\omega)-w_j(r,\omega)=w_j(t-r,\vartheta_r\omega)-w_j(0,\vartheta_r\omega).$$

Consider the process
$$\beta(t)=e^{-\sum_{j=1}^m b_jw_j(t)}.$$
Then, formally, the process $v(t)$ defined by the time change
$$v(t)=\beta(t)\theta(t),$$
satisfies the equation (which depends on a random parameter)
\begin{equation}\frac{dv}{dt}+A_\alpha v+\beta^{-1}u_v\cdot\nabla v=0,\end{equation}
where $u_v$ satisfies (1.3) with $v$ in place of $\theta$.

Then by  similar arguments as the proof of Theorems 3.1-3.3, one can show that for  every $\omega\in \Omega$ the following holds for $s>2(1-\alpha)$:

(i) For all $t_0\in\mathbb{R}$ and $v_0\in H^s$, there exists a unique solution $v\in C([t_0,\infty);H^s)\cap L^2_{\rm{loc}}(t_0,\infty; H^{s+\alpha})$ of equation (6.2) satisfying $v(t_0)=v_0$.

(ii) If such solution is denoted by $v(t,\omega,t_0,v_0)$, the mapping $v_0\rightarrow v(t,\omega;t_0,v_0)$ is continuous in $H^s$ for all $t\geq t_0$.

Then we define
$$\varphi(t,\omega)\theta_0:=\beta(t,\omega)^{-1}v(t,\omega;0,\theta_0), \quad t\geq0.$$
$$S(t,r;\omega)\theta_0:=\beta(t,\omega)^{-1}v(t,\omega;r,\theta_0\beta(r,\omega)), \quad t,r\in \mathbb{R}.$$
\vskip.10in
\th{Theorem 6.1} Fix $\alpha>1/2$. $\varphi(t,\omega)$ is a continuous random dynamical system, and $S(t,r;\omega)$
is a continuous stochastic flow, which is called the stochastic flow associated with the quasi-geostrophic equation driven by multiplicative noise.

\proof By the $\omega$-wise uniqueness of the solution to equation
(6.2)  obtained above, we have that
$$S(t,r;\omega)=S(t,l;\omega)S(l,r;\omega),$$
$$S(t,r;\omega)x = S(t-r,0;\vartheta_r\omega)x,$$
$$\varphi(t+r,\omega)=\varphi(t,\vartheta_r\omega)\circ\varphi(r,\omega),$$
for all $t,l,r\in \mathbb{R}$ and for all $\omega\in\Omega$.
It remains to prove the measurability of $\varphi:\mathbb{R}^+\times \Omega\times H^s\rightarrow H^s$. Since $\varphi(t,\omega)\theta_0=\beta(t,\omega)^{-1}v(t,\omega;0,\theta_0)$, $t\mapsto v(t,\omega;0,\theta_0)$ and $\theta_0\mapsto v(t,\omega;0,\theta_0)$
is continuous, we only need to prove the measurability of $\omega\mapsto v(t,\omega;0,\theta_0)$. By the $\omega$-wise uniqueness of the solutions to (6.2) each subsequence of the convolution approximation $v^n(t,\omega;0,\theta_0)$ we used in the proof of existence of solutions to (6.2) has a subsequence converging to $v(t,\omega;0,\theta_0)$ in $L^2([t_1,t_2],H^s)$ for some $t_1\leq t\leq t_2$. Thus we obtain that the whole sequence of $v^n(t,\omega;0,\theta_0)$  converges to $v(t,\omega;0,\theta_0)$ in $L^2([t_1,t_2],H^s)$, which implies the measurability of  $\omega\mapsto v(t,\omega;0,\theta_0)$.
$\hfill\Box$
\vskip.10in

Fix $t_0<-3$. Now we start with some useful  estimates which lead to the proof of the existence of an absorbing set for the solutions in the space $H^s$ for $s>2(1-\alpha)$. For $s<1$, we choose $p$ such that $\frac{2}{2\alpha-1}<p\leq\frac{2}{1-s}$ and for $s\geq 1$ we take any $p$ satisfying $\frac{2}{2\alpha-1}<p<\infty$. In the following we fix such $p$, we have $H^s\subset L^p$.
\vskip.10in
\th{Lemma 6.2} $$|v(t)|^2\leq |v(t_0)|^2e^{-2\lambda_1 (t-t_0)}, t\geq t_0.$$

Furthermore,
\begin{equation}|v(t+1)|^2+\kappa \int_t^{t+1}|\Lambda^\alpha v|^2dr\leq |v(t_0)|^2e^{-2\lambda_1 (t-t_0)},t\geq t_0.\end{equation}
\proof By above we have $v\in C([t_0,+\infty);H^s)$. Multiplying (6.2) with $v$ and taking the inner product in $L^2$, we have
$$\frac{d}{dt}|v|^2+2\kappa |\Lambda^\alpha v|^2\leq 0.$$
Then Gronwall's lemma yields that
$$|v(t)|^2\leq |v(t_0)|^2e^{-2\lambda_1 (t-t_0)}, t\geq t_0,$$
which implies that
$$|v(t+1)|^2+2\kappa \int_t^{t+1}|\Lambda^\alpha v|^2dr\leq |v(t)|^2\leq |v(t_0)|^2e^{-2\lambda_1 (t-t_0)},t\geq t_0.$$
$\hfill\Box$
\vskip.10in

\th{Lemma 6.3} For $p$ as above, we have\begin{equation}\|v(t)\|_{L^p}\leq \|v(t_0)\|_{L^p}\exp\{-\frac{2\lambda_1}{p}(t-t_0)\}, t\geq t_0.\end{equation}
\proof Multiplying (6.2) with $p|v|^{p-2}v$, taking the inner product in $L^2$ and using  Lemma A.1 we have
$$\frac{d}{dt}\|v\|_{L^p}^p+2\lambda_1 \| v\|_{L^p}^p\leq 0.$$
By Gronwall's lemma, we obtain
(6.4). We can choose a similar approximation $v_n$ as in the proof of Theorem 3.1 to make it rigorously.
$\hfill\Box$
\vskip.10in

\vskip.10in
\th{Lemma 6.4} There exists random radius $r_1(\omega)>0, c_1(\omega)>0,$ and $c_2(\omega)>0$, such that for all $\rho>0$ there exists $t(\omega)\leq -3$ such that the following holds $P$-a.s. : For all $t_0\leq t(\omega)$ and all $\theta_0\in H^s$ with $|\Lambda^s \theta_0|\leq \rho$, the solution $v(t,\omega;t_0,\beta(t_0,\omega)\theta_0)$ with $v(t_0)=\beta(t_0)\theta_0$ satisfies the following inequalities:
\begin{equation}|\Lambda^sv(-1,\omega;t_0,\beta(t_0,\omega)\theta_0)|^2\leq r_1^2(\omega).\end{equation}
\begin{equation}|\Lambda^sv(t,\omega;t_0,\beta(t_0,\omega)\theta_0)|^2\leq c_1(\omega), t\in[-1,0].\end{equation}
\begin{equation}\int_{-1}^0|\Lambda^{s+\alpha}v(t,\omega;t_0,\beta(t_0,\omega)\theta_0)|^2dt\leq c_2(\omega).\end{equation}
\proof To prove Lemma 6.4, first we give the $H^s$-norm estimates of the solutions to (6.2).

[\textbf{$H^s$-norm estimates}] Since the solution of (6.2)  $v\in C([t_0,\infty);H^s)\cap L^2_{\rm{loc}}(t_0,\infty; H^{s+\alpha})$, we obtain for $s_0\leq s$ the following estimate by taking the inner product in $L^2$ with $\Lambda^{s_0}e_k$ for (6.2), multiplying both sides by $\langle v, \Lambda^{s_0} e_k\rangle$, and summing up over $k$:
\begin{equation}\aligned\frac{d}{dt}|\Lambda^{s_0}v|^2+\kappa |\Lambda^{{s_0}+\alpha} v|^2\leq & C\beta^{-1}|\Lambda^{{s_0}+\alpha} v||\Lambda^{{s_0}-\alpha+1+\sigma_1}v|\|v\|_{L^p}\\\leq & \frac{\kappa}{2} |\Lambda^{{s_0}+\alpha} v|^2+C(\beta^{-1}\|v\|_{L^p})^r|\Lambda^{s_0} v|^2\\\leq &\frac{\kappa}{2} |\Lambda^{{s_0}+\alpha} v|^2+C(\|\theta(t_0)\|_{L^p}\beta(t)^{-1}\beta(t_0)\exp\{-\frac{ 2\lambda_1}{p}(t-t_0)\})^r|\Lambda^{s_0} v|^2,\endaligned\end{equation}
where $0<\sigma_1=\frac{2}{p}<2\alpha-1, r:=\frac{2\alpha}{2\alpha-1-\sigma_1}$ as in (3.14). We  used Lemmas 2.7, 2.8 in the first inequality, the interpolation inequality and Young's inequality in the second inequality and  (6.4) in the last inequality. Here the calculation is similar as (3.14) and we omit the details.

By Gronwall's lemma we have for $l\in[-2,-1]$
\begin{equation}\aligned|\Lambda^{s_0}v(-1)|^2\leq & |\Lambda^{s_0}v(l)|^2\exp\{\int_{l}^{-1}C(\|\theta(t_0)\|_{L^p}\beta(\tau)^{-1}\beta(t_0)\exp\{-\frac{ 2\lambda_1}{p}(\tau-t_0)\})^rd\tau\}\\\leq&
|\Lambda^{s_0}v(l)|^2C\exp\{\int_{-2}^{-1}C(\beta(\tau)^{-r}\exp\{-\frac{2r \lambda_1}{p}\tau\})d\tau\|\theta(t_0)\|_{L^p}^r\beta(t_0)^r\exp\{\frac{2r\lambda_1}{p}t_0\}\} .\endaligned\end{equation}
Integrating $l$ over $[-2,-1]$, we obtain
\begin{equation}\aligned|\Lambda^{s_0}v(-1)|^2\leq&
C\int_{-2}^{-1}|\Lambda^{s_0}v(l)|^2dl\exp\{\int_{-2}^{-1}C(\beta(\tau)^{-r}\exp\{-\frac{2r \lambda_1}{p}\tau\})d\tau\|\theta(t_0)\|_{L^p}^r\beta(t_0)^r\exp\{\frac{2r\lambda_1}{p}t_0\}\} .\endaligned\end{equation}
Thus for $s_0= \alpha$, (6.3) yields that
\begin{equation}\aligned|\Lambda^{\alpha}v(-1)|^2\leq&
C|v(t_0)|^2e^{2\lambda_1t_0}\exp\{\int_{-2}^{-1}C(\beta(\tau)^{-r}\exp\{-\frac{2r\lambda_1}{p}\tau\})d\tau\|\theta(t_0)\|_{L^p}^r\beta(t_0)^r\exp
\{\frac{2r\lambda_1}{p}t_0\}\} .\endaligned\end{equation}
By a similar calculation, we also get
\begin{equation}\aligned|\Lambda^\alpha v(-2)|^2\leq&
C|v(t_0)|^2e^{2\lambda_1t_0}\exp\{\int_{-3}^{-2}C(\beta(\tau)^{-r}\exp\{-\frac{2r\lambda_1}{p}\tau\})d\tau\|\theta(t_0)\|_{L^p}^r\beta(t_0)^r\exp
\{\frac{2r\lambda_1}{p}t_0\}\} .\endaligned\end{equation}
Hence by (6.8) and Gronwall's lemma, we have for $t\in [-2,-1]$,
\begin{equation}\aligned|\Lambda^\alpha v(t)|^2\leq& |\Lambda^\alpha v(-2)|^2\exp\{\int_{-2}^{t}C(\|\theta(t_0)\|_{L^p}\beta(\tau)^{-1}\beta(t_0)\exp\{-\frac{ 2\lambda_1}{p}(\tau-t_0)\})^rd\tau\}\\ \leq
&C|v(t_0)|^2e^{2\lambda_1t_0}\exp\{\int_{-3}^{-1}C(\beta(\tau)^{-r}\exp\{-\frac{2 r \lambda_1}{p}\tau\})d\tau\|\theta(t_0)\|_{L^p}^r\beta(t_0)^r\exp\{\frac{2r\lambda_1}{p}t_0\}\} .\endaligned\end{equation}
Moreover, by (6.8), (6.12) and (6.13) we obtain
$$\aligned\int_{-2}^{-1}|\Lambda^{2\alpha} v(l)|^2dl\leq& C|\Lambda^\alpha v(-2)|^2+C\int_{-2}^{-1}(\|\theta(t_0)\|_{L^p}\beta(\tau)^{-1}\beta(t_0)\exp\{-\frac{ 2\lambda_1}{p}(\tau-t_0)\})^rd\tau\\&\sup_{-2\leq t\leq-1}|\Lambda^\alpha v(t)|^2\\\leq &C|\Lambda^\alpha v(-2)|^2+ C|v(t_0)|^2e^{2\lambda_1t_0}\exp\{\int_{-3}^{-1}C(\beta(\tau)^{-r}\exp\{-\frac{2 r \lambda_1}{p}\tau\})d\tau\\&\|\theta(t_0)\|_{L^p}^r\beta(t_0)^r\exp\{\frac{2r\lambda_1}{p}t_0\}\}\\\leq &C|v(t_0)|^2
e^{2\lambda_1t_0}\exp\{\int_{-3}^{-1}C(\beta(\tau)^{-r}\exp\{-\frac{2r \lambda_1}{p}\tau\})d\tau\|\theta(t_0)\|_{L^p}^r\beta(t_0)^r\exp\{\frac{2r\lambda_1}{p}t_0\}\} .\endaligned$$
Therefore by the same arguments as above and a boot-strapping argument, we get for $s>2(1-\alpha)$,
\begin{equation}\aligned|\Lambda^sv(-1)|^2\leq & C|v(t_0)|^2e^{2\lambda_1t_0}\exp\{\int_{-3}^{-1}C(\beta(\tau)^{-r}\exp\{-\frac{ 2r \lambda_1}{p}\tau\})d\tau\|\theta(t_0)\|_{L^p}^r\beta(t_0)^r\exp\{\frac{2r\lambda_1}{p}t_0\}\} .\endaligned\end{equation}

Then by (6.8) we have for $t\in[-1,0]$, $s>2(1-\alpha)$
\begin{equation}\aligned|\Lambda^sv(t)|^2\leq & |\Lambda^sv(-1)|^2\exp\{\int_{-1}^{0}C(\|\theta(t_0)\|_{L^p}\beta(\tau)^{-1}\beta(t_0)\exp\{-\frac{2\lambda_1}{p}(\tau-t_0)\})^rd\tau\}\\\leq&
C|v(t_0)|^2e^{2\lambda_1t_0}\exp\{\int_{-3}^{0}C(\beta(\tau)^{-r}\exp\{-\frac{2r \lambda_1}{p}\tau\})d\tau\|\theta(t_0)\|_{L^p}^r\beta(t_0)^r\exp\{\frac{2r\lambda_1}{p}t_0\}\},\endaligned\end{equation}
and
\begin{equation}\aligned\int_{-1}^{0}|\Lambda^{s+\alpha} v(l)|^2dl\leq &C|\Lambda^s v(-1)|^2+ C\int_{-1}^0(\|\theta(t_0)\|_{L^p}\beta(t)^{-1}\beta(t_0)\exp\{-\frac{2\lambda_1}{p}(t-t_0)\})^rdt\sup_{-1\leq t\leq0}|\Lambda^s v|^2\\\leq&C|v(t_0)|^2e^{2\lambda_1t_0}\exp\{\int_{-3}^{0}C(\beta(\tau)^{-r}\exp\{-\frac{ 2r \lambda_1}{p}\tau\})d\tau\|\theta(t_0)\|_{L^p}^r\beta(t_0)^r\exp\{\frac{2r\lambda_1}{p}t_0\}\} .\endaligned\end{equation}

[\textbf{Absorption in $H^s$ at time $t = -1$}]

 Since $$\lim_{t\rightarrow-\infty}\frac{1}{t}\sum_{j=1}^mb_jw_j(t)=0 \qquad P-a.s.,$$ we have that
$$\beta(t_0)^r\exp\{\frac{2r\lambda_1}{p}t_0\}\rightarrow0 \textrm{ as } t_0\rightarrow-\infty \qquad P-a.s..$$
Then for $|\Lambda^s \theta_0|\leq \rho$, choose $t(\omega)$ such that
$$\|\theta(t_0)\|_{L^p}^r\beta(t_0)^r\exp\{\frac{2r\lambda_1}{p}t_0\}\leq 1,$$
$$|v(t_0)|^2e^{\lambda_1t_0}\leq 1,$$
for all $t_0\leq t(\omega)$.
Hence by (6.14) we get (6.5). (6.6) and (6.7) can be obtained similarly by (6.15) and (6.16).
$\hfill\Box$
\vskip.10in

\th{Lemma 6.5} There exists a random radius $r_2(\omega)>0$, such that for all $\rho>0$ there exists $t(\omega)\leq -1$ such that the following holds $P$-a.s.: For all $t_0\leq t(\omega)$ and all $\theta_0\in H^s$ with $|\Lambda^s \theta_0|\leq \rho$, the solution $v(t,\omega;t_0,\beta(t_0,\omega)\theta_0)$ with $v(t_0)=\beta(t_0)\theta_0$ satisfies the inequality
$$|\Lambda^{s+\alpha}\theta(0,\omega;t_0,\theta_0)|^2\leq r_2^2(\omega).$$
\proof By (6.7) we have for almost every $l\in [-1,0]$, $v(l)\in H^{s+\alpha}$. Then by a similar argument as in the proof of Theorem 3.1 we obtain the solution $v\in L^\infty_{\rm{loc}}([l,\infty); H^{s+\alpha})\cap L^2_{\rm{loc}}([l,\infty);H^{s+2\alpha})$. By a similar estimate as (6.8) we get that
$$\aligned\frac{d}{dt}|\Lambda^{s+\alpha}v|^2+\kappa |\Lambda^{s+2\alpha} v|^2\leq & C\beta^{-1}|\Lambda^{s+2\alpha} v||\Lambda^{s+1+\sigma_1}v|\|v\|_{L^p}\\\leq & \frac{\kappa}{2} |\Lambda^{s+2\alpha} v|^2+C(\beta^{-1}\|v\|_{L^p})^{r}|\Lambda^{s+\alpha} v|^2\\\leq &\frac{\kappa}{2} |\Lambda^{s+2\alpha} v|^2+C(\|\theta(t_0)\|_{L^p}\beta^{-1}(t)\beta(t_0)\exp\{-\frac{2\lambda_1}{p}(t-t_0)\})^{r}|\Lambda^{s+\alpha} v|^2,\endaligned$$
where $\sigma_1,r,p$ are as in (3.14) and we used Lemmas 2.7, 2.8, the interpolation inequality and Young's inequality in the second inequality and Lemma 6.3 in the last inequality.
Therefore Gronwall's lemma implies that
$$\aligned|\Lambda^{s+\alpha}v(0)|^2\leq &|\Lambda^{s+\alpha}v(l)|^2\exp\{\int_l^{0}C(\|\theta(t_0)\|_{L^p}\beta^{-1}(\tau)\beta(t_0)\exp\{-\frac{2\lambda_1}{p}(\tau-t_0)\})^rd\tau\}
\\\leq&|\Lambda^{s+\alpha}v(l)|^2\exp\{\int_{-1}^{0}C(\beta(\tau)^{-r}\exp\{-\frac{ 2r \lambda_1}{p}\tau\})d\tau \|\theta(t_0)\|_{L^p}^r\beta(t_0)^r\exp\{\frac{2r\lambda_1}{p}t_0\}\}.\endaligned$$
Integrating $l$ over $[-1,0]$ and by (6.16) we have
$$\aligned&|\Lambda^{s+\alpha}\theta(0)|^2=|\Lambda^{s+\alpha}v(0)|^2\\\leq & \int_{-1}^0|\Lambda^{s+\alpha}v(l)|^2dl\exp\{\int_{-1}^{0}C(\beta(\tau)^{-r}\exp\{-\frac{ 2r \lambda_1}{p}\tau\})d\tau \|\theta(t_0)\|_{L^p}^r\beta(t_0)^r\exp\{\frac{2r\lambda_1}{p}t_0\}\}\\\leq&C|v(t_0)|^2e^{2\lambda_1t_0}\exp\{\int_{-3}^{0}C(\beta(\tau)^{-r}\exp\{-\frac{2r \lambda_1}{p}\tau\})d\tau\|\theta(t_0)\|_{L^p}^r\beta(t_0)^r\exp\{\frac{2r\lambda_1}{p}t_0\}\} .\endaligned$$
From this and a similar argument as in  the last step of the proof of Lemma 6.4 we have  the absorption of $\varphi$ in $H^{s+\alpha}$ at time $t =0$.$\hfill\Box$
\vskip.10in

Thus by Proposition 2.4 and [8, Corollary 4.6] we obtain the following results.
\vskip.10in
\th{Theorem 6.6} Fix $\alpha>1/2$.  The stochastic flow associated with the quasi-geostrophic equation driven by multiplicative noise (6.1) has a compact stochastic attractor in $H^s$.

Moreover, the Markov semigroup induced by the flow on $H^s$ has an invariant measure $\rho$.
\vskip.10in
\no \textbf{Appendix} In the appendix we will collect some useful results we proved in \cite{RZZ12} for the reader's convenience.

 \th{Lemma A.1}( [17, Lemma 7.4.1] ) For $\alpha\in (0,1)$, and $\theta\in H^1$ with $\Lambda^{2\alpha}\theta\in L^2$, for some $2<p<\infty$, then
$$\int|\theta|^{p-2}\theta(\kappa\Lambda^{2\alpha}-\frac{2\lambda_1}{p})\theta\geq0.$$

\vskip.10in
\th{Proposition A.2} ([17, Proposition 7.4.2]) Let $\alpha>\frac{1}{2}$. Suppose (E.1) holds with $s=1$. Then for $\theta_0\in L^p$, let $\theta$ denote the solution of equation (3.1) with the initial value $\theta_0$ at time $t_0$. Then for $2<p<\infty$, $t>t_0$
$$\aligned&\|\theta(t)\|_{L^{p} }^{p} +\lambda_1\int_{t_0}^{t}\int_{\mathbb{T}^2}|\theta(l)|^{p} d\xi dl\\\leq&\|\theta_0\|_{L^{p} }^{p}
+C_S^p[\frac{1}{2}p(p-1)]^{p/2}\lambda_1^{-\frac{p-2}{2}}\mathcal{E}_0^{p/2}(t-t_0)
+p\int_{t_0}^{t}
\int_{\mathbb{T}^2}|\theta(l)|^{p-2}\theta(l) d\xi dW(l),\endaligned$$
and
$$E\|\theta(t)\|_{L^p}^p\leq \|\theta_0\|_{L^p}^pe^{-\lambda_1(t-t_0)}+\frac{C}{\lambda_1}(1-e^{-\lambda_1(t-t_0)}).$$
\vskip.10in
The following Lemma is a technical result from [12, Lemma 5]. Let $\{X_n\}$ be a sequence of real random variables indexed by $n$. Let $f:\mathbb{Z}^+\rightarrow\mathbb{R}^+$. Define the random variable $T_{\rm{bound}}(\{X_n\},f)$ to be the smallest positive integer such that $m>T_{\rm{bound}}(\{X_n\},f)\Rightarrow |X_m|<f(m).$
\vskip.10in
\th{Lemma A.3} Assume that
$$P(|X_n|\geq \varepsilon n^\delta)\leq \frac{E|X_n|^p}{n^{p\delta}\varepsilon^p}\leq \frac{C}{n^{p\delta-r}\varepsilon^p},$$
for some $\varepsilon,\delta,p,C>0$ and $r\geq 0$. Then $E[T_{\rm{bound}}(\{X_n\},\varepsilon\delta^n)]^q<\infty$ for $q\in (0,p\delta-(1+r))$.

\vskip.10in
\th{Acknowledgement.} The authors would like to thank Professor M. R${\ddot{o}}$ckner for valuable discussions and suggestions.

\end{document}